\documentclass[a4paper,10pt,leqno,twoside]{article}

\usepackage[english]{babel}
\usepackage{amsmath, amssymb , latexsym, amsthm, epic, epsfig, rotating, fancyheadings, pifont, empheq}

\newcommand{\alphlist}{\begin{list}{(\alph{enumi})}{\usecounter{enumi}}}
\newcommand{\romanlist}{\begin{list}{(\roman{enumi})}{\usecounter{enumi}}}
\newcommand{\listend}{\end{list}}

\newcommand{\ld}{\ensuremath{,\ldots,}}
\newcommand{\ssq}{\ensuremath{\subseteq}}
\newcommand{\smin}{\ensuremath{\setminus}}

\newcommand{\ra}{\ensuremath{\rightarrow}}

\newcommand{\follows}{\ensuremath{\Rightarrow}}

\newcommand{\N}{\ensuremath{\mathbb{N}}} 
\newcommand{\R}{\ensuremath{\mathbb{R}}}
\newcommand{\Z}{\ensuremath{\mathbb{Z}}}
\newcommand{\Q}{\ensuremath{\mathbb{Q}}}

\newcommand{\kreis}{\ensuremath{\mathbb{T}^{1}}}
\newcommand{\ntorus}[1][2]{\ensuremath{\mathbb{T}^{#1}}}

\newcommand{\nLim}{\ensuremath{\lim_{n\rightarrow\infty}}}

\newcommand{\xLim}{\ensuremath{\lim_{x\rightarrow\infty}}}

\newcommand{\nfolge}[1]{\ensuremath{(#1)_{n\in\mathbb{N}}}}

\newcommand{\ifolge}[1]{\ensuremath{(#1)_{i\in\mathbb{N}}}}

\newcommand{\supp}{\ensuremath{\mathrm{supp}}}

\newcommand{\ncap}{\ensuremath{\bigcap_{n\in\N}}}

\newcommand{\ktel}{\ensuremath{\frac{1}{k}}}

\newcommand{\ntel}{\ensuremath{\frac{1}{n}}}

\newcommand{\ptel}{\ensuremath{\frac{1}{p}}}

\newcommand{\halb}{\ensuremath{\frac{1}{2}}}

\newcommand{\viertel}{\ensuremath{\frac{1}{4}}}

\newcommand{\thx}{\ensuremath{(\theta,x)}}
\newcommand{\thom}{\ensuremath{\theta + \omega}}

\newcommand{\fth}{\ensuremath{f_{\theta}}}
\newcommand{\fthx}{\ensuremath{f_{\theta}(x)}}

\title{\textsc{Strange non-chaotic attractors in quasiperiodically forced
  circle maps}} \author{Tobias J\"ager}
\newtheorem{definition}{Definition}[section]
\newtheorem{thm}[definition]{Theorem}

\newtheorem{lem}[definition]{Lemma}
\newtheorem{cor}[definition]{Corollary}  
\newtheorem{prop}[definition]{Proposition}

\theoremstyle{plain} \newtheorem{claim}[definition]{Claim} 

\newtheorem{bem}[definition]{Remark}

\numberwithin{equation}{section}
\numberwithin{figure}{section}

\newcommand{\nofolge}[1]{\ensuremath{(#1)_{n\in\mathbb{N}_0}}}
\newcommand{\eps}{\ensuremath{\varepsilon}}
\newcommand{\Leb}{\ensuremath{\mathrm{Leb}}}

\newcommand{\sltr}{\ensuremath{\textrm{SL}(2,\mathbb{R})}}
\newcommand{\Id}{\ensuremath{\textrm{Id}}}

\begin{document}
\setlength{\evensidemargin}{0.06\textwidth}  

\pagestyle{fancy}

\chead[Tobias J\"ager]{\textsc{SNA in
    quasiperiodically forced circle maps}}
\lhead[\arabic{page}]{}
\rhead[]{\arabic{page}}

\cfoot{}

\maketitle

\begin{abstract}
The occurrence of strange non-chaotic attractors (SNA) in quasiperiodically
forced systems has attracted considerable interest over the last two decades,
in particular since it provides a rich class of examples for the possibility
of complicated dynamics in the absence of chaos. Their existence was
discovered in the early 1980's, independently by Herman \cite{herman:1983} for
quasiperiodic \sltr-cocycles and by Grebogi {\em et al}
\cite{grebogi/ott/pelikan/yorke:1984} for so-called {\em `pinched skew
products'}. However, except for these two particular classes there are still
hardly any rigorous results on the topic, despite a large number of numerical
studies which all confirmed the widespread existence of SNA in
quasiperiodically forced systems.

Here, we prove the existence of SNA in quasiperiodically forced circle maps
under rather general conditions, which can be stated in terms of ${\cal
C}^1$-estimates. As a consequence, we obtain the existence of strange
non-chaotic attractors for parameter sets of positive measure in suitable
parameter families.  Further, we show that the considered systems have minimal
dynamics. The results apply in particular to a forced version of the Arnold
circle map. For this particular example, we also describe how the first Arnold
tongue collapses and looses its regularity due to the presence of strange
non-chaotic attractors and a related unbounded mean motion property.
\end{abstract}


\section{Introduction}

In 1984, Grebogi {\em et al} \cite{grebogi/ott/pelikan/yorke:1984} introduced
a class of quasiperiodically forced (qpf) interval maps which exhibit
non-continuous invariant graphs with negative (vertical) Lyapunov
exponents. As these objects attract a set of initial conditions of positive
measure and combine a complicated structure with non-chaotic dynamics (in
particular zero topological entropy), they are commonly referred to as {\em
strange non-chaotic attractors} (SNA).  Already one year earlier, Herman
\cite{herman:1983} had proved the existence of such SNA in certain parameter
families of qpf circle diffeomorphisms that are induced by the projective
action of $\sltr$-cocycles over an irrational rotation (see also
\cite{haro/puig:2006}).

In the following years, the phenomenom attracted a considerable amount of
interest, and a large number of numerical studies indicated that the existence
of SNA is quite common in quasiperiodically forced systems
(\cite{prasad/negi/ramaswamy:2001} gives a good overview and further
reference). However, despite all efforts rigorous results remained rare, and
in particular the two classes of examples mentioned above remained the only
ones for which the existence of SNA could be proved rigorously. Only recently
some further progress was made, as the author described the creation of SNA in
non-smooth bifurcations of invariant curves, which take place in qpf interval
maps \cite{jaeger:2006b} (but only at isolated parameter values).
\medskip

The aim of this article is two-fold. First, we show that once the skew-product
structure is given, which is usually motivated by the physical
context of the model, the existence of SNA in qpf circle maps is a phenomenom
which is both `robust' and `non-degenerate'. To make this more precise, we denote
by $\textrm{Diff}_0(\ntorus)$ the set of all diffeomorphisms of the two-torus
which are homotopic to the identity and by $\pi_i$ the projection to the
respective coordinate. Further, for any $\omega \in \kreis$ we let
$R_\omega\thx = (\thom,x)$. Then, as a consequence of our results, we obtain the
following: 
\begin{quote} \em
  Let ${\cal F} := \{ F \in \mathrm{Diff}_0(\ntorus) \mid \pi_1 \circ F =
  \pi_1 \}$. Then there exists a non-empty set ${\cal U} \ssq {\cal F}$, which
  is ${\cal C}^1$-open in ${\cal F}$ and has the following property:

  For any $F\in{\cal U}$ there exists a set $\Omega_F \ssq \kreis$ of positive
  Lebesgue measure, such that for any $\omega \in \Omega_F$ the map
  $f=R_\omega \circ F$ is minimal and has a strange non-chaotic attractor.
\end{quote}
A more precise characterisation of the set ${\cal U}$
in the above statement, in terms of explicit ${\cal C}^1$-estimates, is provided
by Theorem~\ref{thm:basic} and/or Theorem~\ref{thm:refined} below. 
\medskip

Our second objective is to apply our methods to a particular model, which is
well-known from the literature, namely the qpf Arnold circle
map
\begin{equation}
  \label{eq:arnold} \thx \ \mapsto \ \left(\thom,x+\tau+a\sin(2\pi x) +
  b\cos(2\pi\theta)^d \right) \ .
\end{equation}
Here $\tau\in\kreis,a\in[0,1/2\pi],b\in\R$ and $d$ is an odd positive
integer. This example was proposed by Ding {\em et al}
\cite{ding/grebogi/ott:1989} as a simple model of an oscillator forced at two
or more incommensurate frequencies, and has been intensively studied
numerically since%
\footnote{In the numerical studies usually $d=1$. However, as mentioned in
  \cite{ding/grebogi/ott:1989}, any real-analytic forcing function is of more
  or less equal interest.}
(see, for example,
\cite{feudel/kurths/pikovsky:1995,chastell/glendinning/stark:1995,%
glendinning/feudel/pikovsky/stark:2000,osinga/wiersig/glendinning/feudel:2001%
,stark/feudel/glendinning/pikovsky:2002}).
Provided $d$ is chosen sufficiently large, we show that there exist
rotation numbers $\omega$ for which (\ref{eq:arnold}) exhibits SNA on a set of
positive measure in the $(\tau,a,b)$-parameter space (see
Corollary~\ref{cor:arnold}). 

Particular attention in the study of (\ref{eq:arnold}) has been given to the
structure of the Arnold tongues, which are subsets of the parameter space on
which the rotation number stays constant. In
\cite{glendinning/feudel/pikovsky/stark:2000}, the authors observe that the
Arnold tongue corresponding to rotation number zero seems to collapse in some regions of the
parameter space. In Section~\ref{ArnoldTongues}, we prove that this happens at
least for large $d$. In addition, we show that the boundaries of the zero
tongue do not depend analytically on the parameter $\beta$ in this case.
\medskip

We want to mention that the approach employed here is inspired by the one of
Bjerkl\"ov in \cite{bjerkloev:2005a}. The latter was developed in the setting
of quasiperiodic Schr\"odinger cocycles, but its techniques are basically
non-linear, which allows us to adapt and to apply them to the non-linear
setting. Similar ideas have also been used earlier by Young
\cite{young:1997} to prove positive Lyapunov exponents for certain
quasiperiodic $\sltr$-cocycles.
\medskip

\noindent
\textbf{Acknowledgements.} I would like to thank Kristian
Bjerkl\"ov, for inspiration and stimulating discussions, as well as
J.-C.~Yoccoz and the Coll\`ege de France for their hospitality during a
two-year visit.  This work was supported by a research fellowship of
the German Research Council (DFG).


\subsection{Notation}
\noindent
Let $\kreis := \R / \Z$ and denote by $\pi_i : \ntorus \ra \kreis$ the
projection to the respective coordinate. A \emph{quasiperiodically forced
(qpf) circle homeomorphism/diffeomorphism} is a homeomorphism/diffeomorphism
$f : \ntorus \ra \ntorus$ which is of the form
\begin{equation} \label{eq:skew}
   f: \thx \mapsto (\thom,f_\theta(x))
\end{equation}
where $\omega \in \kreis \smin \Q$ and the \emph{fibre maps} $\fth$ are
defined by $\fth(x) = ´\pi_2 \circ f(\theta,x)$. Derivatives with respect to
$\theta$ or $x$ will be denoted by $\partial_\theta$ and $\partial_x$,
respectively. Further, we use the notation
\[
f^n_\theta(x) \ := \
\pi_2 \circ f^n\thx \ \ \ \ \ \forall n \in \Z \ .
\]
Note that this implies $f^{-1}_\theta = (f_{\theta-\omega})^{-1}$. For any
$a,b\in\kreis$, we denote by 
\[
[a,b] \ := \ \{x \in\kreis\mid a \leq x \leq b\}
\]
the interval of all points $x \in \kreis$ which lie between $a$ and $b$ in the
counterclockwise direction, similarly for open intervals. Note that thus
$[b,a] = \kreis \smin (a,b)$. For two points $x,y\in \kreis$, we denote the
usual Euclidean distance on the circle by $d(x,y)$. We will also
use the notation $y-x$ in order to denote the distance between $x$ and $y$ in
the counterclockwise direction, i.e.\ the length of the interval $[x,y]$.

If $\varphi,\psi : \kreis \ra \kreis$ are two measurable functions, we let
\[
[\varphi,\psi] \ := \ \{\thx \in \ntorus \mid x \in
[\varphi(\theta),\psi(\theta)] \}
\]
 For any initial point $(\theta_0,x_0) \in
\ntorus$ we denote its orbit by $(\theta_k,x_k)_{k\in\Z}$, that is
\[
(\theta_k,x_k) \ := \ f^k(\theta_0,x_0) \ .
\]

\subsection{Some preliminaries}

 An \emph{invariant graph} is a measurable function $\varphi :
\kreis \ra \kreis$ which satisfies
\[
f_\theta(\varphi(\theta)) \ = \ \varphi(\theta+\omega) \ \ \ \ \ \forall
\theta \in \kreis .
\]
This implies that the corresponding point set $\Phi := \{
(\theta,\varphi(\theta)) \mid \theta \in \kreis\}$ is $f$-invariant. 
The Lyapunov exponent of an invariant graph $\varphi$ is defined as
\[
\lambda(\varphi) \ = \ \int_{\kreis} \log |\partial_x f_\theta(\varphi(\theta))|
\ d\theta \ .
\]
We call a non-continuous invariant graph a \emph{strange non-chaotic
  attractor} (SNA) if its Lyapunov exponent is negative and a \emph{strange
  non-chaotic repeller} (SNR) if it is positive.

A convenient criterium for the existence of SNA involves pointwise Lyapunov
exponents, forwards and backwards in time. These are given by
\[
\lambda^+(\theta,x) \ = \ \limsup_{n \ra \infty} \ntel |\log \partial_x
f^n_\theta(x)|
\]
and
\[
\lambda^-(\theta,x) \ = \ \limsup_{n \ra \infty} \ntel |\log \partial_x
f^{-n}_\theta(x)| \ .
\]
A point $(\theta,x) \in \ntorus$ (or more precisely its orbit) which has a
positive Lyapunov exponent both forwards and backwards in time is called a
\emph{sink-source-orbit}. The existence of such orbits implies the existence
of SNAs: 
\begin{prop}[\cite{jaeger:2006b}] \label{prop:sinksourcesna}
  Suppose $f$ is a quasiperiodically forced circle diffeomorphism which has a
  sink-source-orbit. Then $f$ has both a SNA and a SNR.
\end{prop}
The proof in \cite{jaeger:2006b} is given for qpf monotone interval maps, but
using \cite[Theorem~4.1]{furstenberg:1961} it can easily be adapted to qpf
circle diffeomorphisms.
\medskip

The fibred rotation number of a qpf circle homeomorphism is defined as
$\rho(f) = \rho(F) \bmod 1$, where $F : \kreis \times \R \xhookleftarrow{}$ is
a lift of $F$ and
\begin{equation}
  \label{e.rotnum}
\rho(F) \ := \ \nLim \ntel (F^n_\theta(x)-x) \ .
\end{equation}
This limit always exists and is independent of $\thx$
\cite{herman:1983}. Concerning the behaviour of the fibred rotation number
with respect to strictly monotone perturbations, we will make use of the
following: 
\begin{prop}[\cite{bjerkloev/jaeger:2006a}] \label{p.strict-monotonicity}
  Suppose a qpf circle homeomorphism $f$ is minimal. Let $F$ be a lift of
  $f$ and $F_\eps\thx := (\thom,F_\theta(x)+\eps)$. Then the mapping $\eps
  \mapsto \rho(F_\eps)$ is strictly monotone in $\eps = 0$. 
\end{prop}
In fact, the statement given in \cite{bjerkloev/jaeger:2006a} is more general:
The assertion of the proposition is true whenever $f$ has no invariant strip,
which is the appropriate analogue of a periodic orbit in this context (see
\cite{jaeger/stark:2006} or \cite{jaeger/keller:2006} for
the precise definition). Since invariant strips are always compact invariant
strict subsets of \ntorus, the above version follows immediately. 
\medskip

Finally, we will need a result concerning the uniqueness of the minimal set: 
\begin{prop}[\cite{beguin/crovisier/jaeger/leroux:2006b}] \label{p.unique-minimal-set}
  Suppose a qpf circle homeomorphism $f$ is transitive. Then it has a unique
  minimal set. 
\end{prop}


\section{Main results} 

\subsection{The existence of SNA and a first application}
  \label{MainAssumptions}

In the following, we will formulate a number of assumptions which are used in
the statements of our main results. It is important to note that none of them
involves the rotation number $\omega$ on the base, since this will later be
seen as a free parameter of the system. Thus, all the following conditions
should be understood as assumptions on a collection of fibre maps
$(f_\theta)_{\theta\in\kreis}$. Equivalently, the latter might be considered as
a map $F$ which satisfies $\pi_1 \circ F = \textrm{Id}$, as in highlighted
statement in the introduction, such that $F\thx = (\theta,f_\theta(x))$.

\ \\ \emph{I. Regions in the phase space.} Suppose ${\cal I}_0 \ssq \kreis$ is
a finite union of ${\cal N}$ disjoint open intervals $I_0^1 \ld I_0^{\cal
N}$. We will refer to ${\cal I}_0$ as the \emph{first critical
region}. Further, suppose that $E=[e^-,e^+]$ and $C=[c^-,c^+]$ are two non-empty, compact and
disjoint intervals of positive length in \kreis. We will call $E$ the
\emph{expanding} and $C$ the \emph{contracting interval}, motivated by the
bounds on the derivatives given below. The first condition we require is a
strong forward invariance of the contracting interval outside of the critical
region:
\begin{equation}  \label{eq:Cinvariance} \tag{${\cal A}1$}
  f_\theta(\mbox{cl}(\kreis \smin E)) \ \ssq \ \mbox{int}(C) \ \ \ \ \ \forall
  \theta \notin {\cal I}_0 \ .
\end{equation}
Note that this implies
\begin{equation}  \label{eq:Einvariance} \tag{${\cal A}1'$}
  f^{-1}_\theta(\mbox{cl}(\kreis \smin C)) \ \ssq \ \mbox{int}(E) \ \ \ \ \ \forall \theta \notin
  {\cal I}_0+\omega \ .
\end{equation}

\ \\ \emph{II. Bounds on the derivatives.} Let $\alpha =
(\alpha_l,\alpha_c,\alpha_e,\alpha_u) \in \R^4$ satisfy 
\begin{eqnarray*} 
& 0 \ < \ \alpha_l \ < \ \alpha_c \ < \  1 \ < \ \alpha_e \ < \ \alpha_u
\end{eqnarray*}
and suppose the following estimates hold:
\begin{equation} \tag{${\cal A}2$}
 \label{eq:bounds1}  \alpha_l \ < \ \partial_x\fthx \ < \ \alpha_u \hspace{2eM} \quad
  \forall \thx \in \ntorus \ ;
\end{equation}
\begin{equation} \label{eq:bounds2} \hspace{3eM}  \partial_x\fthx \ > \ \alpha_e
  \hspace{4.1eM} \quad \forall \thx \in \kreis \times  \tag{${\cal A}3$}
  E \ ; \end{equation}
\begin{equation} \label{eq:bounds3} \hspace{3eM} \partial_x\fthx \ < \ \alpha_c
  \hspace{4.1eM} \quad \forall \thx \in \kreis \times C \ .    \tag{${\cal A}4$}
 \end{equation}
$\alpha_e$ and $\alpha_c$ will be referred to as the
{\em expansion} and {\em contraction constants}, $\alpha_l$ and $\alpha_u$ as
the {\em lower}
and {\em upper bounds} (on the derivatives $\partial_x\fth$). 

Simply due to compactness, there also exists a global bound for the derivative
w.r.t.\ $\theta$, i.e.\ a constant $S>0$ such that 
\begin{equation} \label{eq:bounddth}  \tag{${\cal A}5$}
   |\partial_\theta \fthx| \ < \ S \ \ \ \ \ \forall \thx \in
  \ntorus \ .
\end{equation}

\ \\ \emph{III. Transversal Intersections.} The last property we will need is
the fact that for each connected component $I^\iota_0$ of the first critical
region ${\cal I}_0$, the set $f(I^\iota_0 \times C)$ crosses the expanding
strip $\kreis \times E$ in a `nice' transversal intersection, either upwards
or downwards. This is ensured by the following: First, we suppose that
\begin{equation} 
   |\partial_\theta \fthx| \ > \ s \ \ \ \ \ \forall \thx \in {\cal I}_0
    \times \kreis \ \label{eq:s} \tag{${\cal A}6$}
\end{equation}
for some constant $s$ with $0 < s < S$.  In particular, this implies that the
sign of $\partial_\theta \fthx$ is constant on every connected component
$I^\iota_0 \times \kreis$ of ${\cal I}_0\times \kreis$.  We speak of an
\emph{upwards crossing} if it is positive and of a {\em downwards crossing} if
it is negative. Secondly, we assume that
\begin{equation} \label{eq:crossing}  \tag{${\cal A}7$}
 \begin{array}{l} \exists!\theta_\iota^1 \in I_0^\iota \textrm{ with  }
   f_{\theta_\iota^1}(c^+) = e^- \textrm{ \ and } \\ \exists! \theta_\iota^2 \in
  I_0^\iota \textrm{ with } f_{\theta_\iota^2}(c^-) = e^+ \ . \end{array}
\end{equation}
This ensures that the image of $I_0^\iota \times C$ crosses the strip
$(I^\iota_0+\omega) \times E$ exactly once and does not `wind around the
torus' several times.  Note that with respect to the canonical ordering inside
the interval $I^\iota_0$, the point $\theta_\iota^1$ lies on the right of
$\theta_\iota^2$ if the crossing is upwards and on the left of
$\theta_\iota^2$ if it is downwards.
\bigskip

Now we can state the first main result. The proof is given in Section
\ref{BasicMechanism}.
\begin{thm} \label{thm:basic}
  Suppose $(f_\theta)_{\theta\in\kreis}$  satisfies
  (\ref{eq:Cinvariance})--(\ref{eq:crossing}). Further assume that
  $$\alpha_c^{-1} = \alpha_e = \alpha^{\frac{2}{p}} \quad \textrm{ and } \quad
  \alpha_l^{-1} = \alpha_u = \alpha^{p}$$ for some $p \in \N$. Let $\eps_0 :=
  \max_{\iota=1}^{\cal N} |I^\iota_0|$ and fix $\delta>0$. Then there exists
  strictly positive constants $c_0=c_0(\delta,p,s,S,{\cal N})$ and
  $\alpha_0=\alpha_0(\delta,p,s,S,{\cal N})$ with the following property:\bigskip

  If $\eps_0 < c_0$ and $\alpha > \alpha_0$, then there exists a set $\Omega
  \ssq \kreis$ of measure
  \[
   \Leb(\Omega) \ \geq \ 1 - \delta \ ,
  \]
  such that for all $\omega \in \Omega$ the system
  \[
   \thx \mapsto (\thom,\fthx)
  \]
  has a sink-source-orbit, and consequently a SNA and a SNR. In addition, the
  dynamics are minimal.
\end{thm}

\begin{bem} \alphlist
\item Since all the conditions of the theorem are ${\cal C}^1$-open in ${\cal
  F}$, the highlighted statement in the introduction is an immediate
  consequence.
\item 
Suppose that a qpf circle diffeomorphism $f$ is minimal and has a SNA, as in
  the assertion of the theorem. Then it also has the property that its {\em
  `deviations from the average rotation'}
\begin{equation}
  \label{e.deviations}
F^n_\theta(x)-x-n\rho(F)
\end{equation}
are unbounded. This follows from a classification result for qpf circle
homeomorphisms, which we want to discuss briefly. 

 If the quantities in (\ref{e.deviations}) are uniformly bounded in $n,\theta$
and $x$, then a direct analogue to Poincar\'e's classification of circle
homeomorphism holds \cite{jaeger/stark:2006}: Either $f$ is semi-conjugate to
an irrational torus translation, or there exists an invariant strip. The
latter replace periodic orbits and are defined as compact invariant sets which
intersect every fibre $\{\theta\}\times\kreis$ in a finite number of intervals
and have certain additional regularity properties (a precise definition is
contained in \cite{jaeger/stark:2006} or
\cite{beguin/crovisier/jaeger/leroux:2006b}).

Since $f$ is minimal it cannot have an invariant strip (such sets are always
strict subsets of the torus), since it has an SNA it cannot be semi-conjugate
to an irrational torus translation (in this case there are no invariant
graphs). Consequently, the two alternatives in the case of bounded deviations
are ruled out, and the quantities in (\ref{e.deviations}) have to be
unbounded. 
\item There exists a mechanism for the creation of SNA which is very similar
  to the one studied here, but which leads to SNA which are the
  semi-continuous boundary graphs of invariant strips. In particular, the
  dynamics are not minimal and the deviations from the constant rotation
  (\ref{e.deviations}) remain bounded. This mechanism is described in
  \cite{bjerkloev:2005} and \cite{jaeger:2006b}.  \listend
\end{bem}

In order to give some explicit examples to which the above theorem applies,
denote by $\gamma:\kreis \ra (-1/2,1/2)$ the lift of the identity map on
\kreis. Then $\pi \circ \gamma = \Id_{\kreis}$, where $\pi : \R \to \kreis$ is
the canonical projection. Further, given any $p\geq 2$ define $a : \R \to \R$
by
\begin{equation} \label{e.ap-def}
a_p(x) \ := \ \int_0^x \frac{1}{1+|x|^p} \ dx \  .
\end{equation}
Of course, for $p=2$ this just yields the arcus tangent. For a given parameter
$\alpha\in\R^+$ and $x \in \kreis$, let
\begin{equation}
 h_\alpha(x) \ := \ \pi\left(\frac{a_p(\alpha \gamma(x))}{2 a_p(\alpha/2)}\right) \ .
\end{equation}
It is easy to check that for all $\alpha$ the map $h_\alpha$ is a
diffeomorphism of the circle. Finally, let $g\in \textrm{Diff(\kreis)}$ be
such that
\begin{eqnarray}
&& g^{-1}(\{1/2\}) \ \textrm{is a finite and non-empty set} \ ; \label{e.g} \\
&&  g'(\theta) \ \neq \ 0 \quad \forall x \in g^{-1}(\{1/2\}) \ . \label{e.g'}
\end{eqnarray}
For example, one could choose $g(\theta) = \beta \cos(2\pi\theta)$ for any
$\beta > \halb$. Then Theorem~\ref{thm:basic} implies the following
\begin{cor}
  \label{cor:firstexample}
Suppose $h_\alpha$ and $g$ are chosen as above and $\delta>0$ is fixed. Then
there exists a constant $\alpha_0 = \alpha_0(\delta,p,g)$ with the following
property:

If $\alpha \geq \alpha_0$, then there exists a set $\Omega\ssq \kreis$ of
measure $\Leb(\Omega) \geq 1-\delta$, such that for any $\omega\in \Omega$ the
system
\begin{equation} \label{e.firstexample}
\thx \ \mapsto \ (\thom,h_\alpha(x)+g(\theta))
\end{equation}
has a sink-source-orbit and consequently a SNA and a SNR. In addition, the
dynamics are minimal.
\end{cor}
The proof is given in Section~\ref{FirstExample}~. 
\begin{bem}
  Let $c_p := \xLim a_p(x)$ and suppose $\tilde{h}_\alpha$ is the map which is
  obtained by projecting the mapping $\bar{\R} \xhookleftarrow{},\ x \mapsto
  \alpha^2x$ to the circle via the change of variables $x \mapsto
  \pi(a_p(x)/2c_p)$. Then the preceding corollary remains true if $h_\alpha$
  is replaced by $\tilde{h}_\alpha$. The proof in Section~\ref{FirstExample}
  can be adapted easily.

  However, in this case the map $\thx \mapsto
  (\thom,\tilde{h}_\alpha(x)+g(\theta))$ is the projective action of the
  $\mbox{SL}(2,\R)$-cocycle $$\kreis \times \R^2 \xhookleftarrow{} \quad
  ,\quad (\theta,v) \mapsto (\thom,A(\theta) v)$$ with $$A(\theta) =
  R_{g(\theta)} \circ \left(\begin{array}{cc}\alpha & 0 \\ 0 & \alpha^{-1}
  \end{array}\right) \ , $$
  where $R_\phi$ denotes the rotation matrix with angle $\phi$. This means
  that, at least in the case of an analytic forcing function $g$ and except
  for the minimality, similar statements can be derived from
  classical results on $\mbox{SL}(2,\R)$-cocycles, for example in
  \cite{herman:1983}. This is not true for the parameter family
  (\ref{e.firstexample}).
\end{bem}


\subsection{A refined result for the quasiperiodically forced Arnold circle map}

\label{Refined}

The statement of Theorem \ref{thm:basic} can be circumscribed by saying that
SNA occur whenever the fibre maps are `sufficiently hyperbolic', meaning that
the expansion and contraction constants provided by (\ref{eq:bounds2}) and
(\ref{eq:bounds3}) are large enough. However, concerning the forced Arnold
circle map (\ref{eq:arnold}), this constitutes a problem. In the realm of
invertibility, meaning for $a \leq 1/2\pi$, the derivative of the fibre maps
is always bounded by 2. For the contraction, the situation is similar: While
the derivative at $x=\halb$ goes to zero as $a$ goes to one, a strong
contraction only takes place on a very small neighbourhood of the point $\halb$. For any interval of
fixed length, the uniform contraction rate will always remain bounded.

In order to overcome this obstruction and to obtain a result which applies to
the qpf Arnold circle map, we have to make use of additional information on
the forcing function $\theta \mapsto \cos(2\pi\theta)^d$, namely of the fact
that for large $d$ its derivative almost vanishes on a large part of the phase
space. This is done via the following assumption.
\medskip

Suppose ${\cal I}_0'\ssq \kreis$ is the disjoint union of at most ${\cal N}$
open intervals and let $s' \in (0,S)$. Then assume that
\begin{equation} \label{eq:refinedbounddth}
  \tag{${\cal A}8$} {\cal I}_0 \ssq {\cal I}_0' \quad \textrm{and} \quad
  |\partial_\theta \fthx| \ < \ s' \ \ \forall \thx \in (\kreis \smin
  {\cal I}_0') \times C \ .
\end{equation}
The refined version of Theorem~\ref{thm:basic} now reads as follows:

 \begin{thm} \label{thm:refined}
  Suppose $(f_\theta)_{\theta\in\kreis}$ satisfies
  (\ref{eq:Cinvariance})--(\ref{eq:refinedbounddth}) and
  $$\alpha_c^{-1} = \alpha_e = \alpha^{\frac{2}{p}} \quad \textrm{ and } \quad
  \alpha_l^{-1} = \alpha_u = \alpha^{p}$$ for some $p \in \N$. Let $\eps_0 :=
  \max_{\iota=1}^{\cal N} |I^\iota_0|$ and fix $\delta>0$. Further, assume
  there exist constants $A,d > 0$ such that
  \begin{eqnarray}
    S & < & A\cdot d \ , \label{e.S<Ad}\\ s & > & \sqrt{d}/A \ , \label{e.s>d/A} \\
    \eps_0 & < &A/\sqrt[3]{d} \ .\label{e.eps<1/Ad}
  \end{eqnarray}

  Then there exist strictly positive constants $c_0=c_0(\delta,\alpha,p,{\cal
  N})$ and $d_0=d_0(\delta,\alpha,p,{\cal N},A)$ with the following
  property:\bigskip

  If $\frac{s'}{s} < c_0$ and $d \geq d_0$, then there exists a set $\Omega
  \ssq \kreis$ of measure
  \[
   \Leb(\Omega) \ \geq \ 1 - \delta \ ,
  \]
  such that for all $\omega \in \Omega$ the system
  \[
  \thx \mapsto (\thom,\fthx)
  \]
  has a sink-source-orbit and consequently a SNA and a SNR. In addition, the
  dynamics are minimal.
\end{thm}

Now suppose $h$ is an orientation-preserving diffeomorphism of the circle,
such that there exists disjoint closed intervals $C,E\ssq \kreis$ which
satisfy 
\begin{equation} \label{e.h-cond1}
\sup_{x\in C} h'(x) \ < \  1 \quad , \quad \inf_{x\in E} h'(x) \ > \ 1 \quad
\end{equation}
and
\begin{equation} \label{e.h-cond2}
h(\mbox{cl}(E^c)) \ \ssq \ \mbox{int}(C) \ .
\end{equation}
For example, this holds whenever $h$ has exactly two fixed points and exactly
two points of inflexion.
\begin{cor} \label{cor:refined}
  Suppose $h$ satisfies (\ref{e.h-cond1}) and (\ref{e.h-cond2}) and $\delta>0$
  is fixed. Then there exist constants $d_0 = d_0(\delta,h)$ and
  $\eps=\eps(\delta,h)$ with the following property:

  If $d \geq d_0$ and $b \in [1-\eps,1+\eps]$, then there exists a
  set $\Omega \ssq \kreis$ of measure $\Leb(\Omega) \geq 1-\delta$, such that
  for any $\omega \in \Omega$ the system
  \begin{equation} \label{e.cor-refined}
       \thx \ \mapsto \ (\thom,h(x)+b\cos(2\pi\theta)^d)
  \end{equation}
  is minimal and has a SNA and a SNR.
\end{cor}
The proof is given in Section~\ref{Proof-Cor-Refined}~.
\begin{bem}
\alphlist 
\item Corollary~\ref{cor:refined} applies in particular to $h(x) = x + \tau +
a\sin(2\pi x)$ whenever $0\leq \tau < a < 1/2\pi$. Thus, we obtain the
existence of SNA for the qpf Arnold circle map (\ref{eq:arnold}). We
reformulate the result in Corollary~\ref{cor:arnold} below.
\item We remark that the above statement remains true if $\cos(2\pi\theta)^d$
  is replaced by other forcing functions depending on a parameter $d$, as long
  as these show a similar scaling behaviour. For example, one could take
  $g_d(\theta) = \left(\frac{1+\sin(2\pi\theta)}{2}\right)^d$. In this case
  $d\in\R^+$ can be chosen either very large or very small in order to ensure
  the existence of SNA. The proof of the corollary in
  Section~\ref{Proof-Cor-Refined} can be adapted accordingly. However, the
  symmetry $\cos(2\pi(\theta+\halb))^d = -\cos(2\pi\theta)^d$ will play an important
  role in Section~\ref{ArnoldTongues}, such that we concentrate on this choice
  of the forcing function.  \listend
\end{bem}

 In the literature, a typical point of view is to consider $\omega$ and $d$ as
fixed and to view (\ref{eq:arnold}) as a three-parameter family depending on
$\tau,a$ and $b$. As a simple consequence of Fubini's Theorem we obtain
\begin{cor} \label{cor:arnold} There is a constant $d_0 > 0$, such that
  for any $d\geq d_0$ there exists a set of positive measure $\Omega \ssq
  \kreis$ with the following property:

  For each $\omega \in \Omega$ there exists a set of positive measure
  $B_\omega \ssq \kreis \times [0,1/2\pi]\times \R$, such that for all
  $(\tau,a,b)\in B_\omega$ the qpf Arnold circle map (\ref{eq:arnold}) is
  minimal and has a SNA and a SNR.
\end{cor}
Of course, similar statements hold if one likes to consider (\ref{eq:arnold})
as parameter family only depending on one or two parameters, while the
other(s) are fixed.


\subsection{Collapsing of the first Arnold tongue} \label{ArnoldTongues}

In this section, we explain the consequences of our results for the structure
of the first Arnold tongue. We denote the qpf Arnold circle map
(\ref{eq:arnold}) with parameters $\tau,a$ and $b$ by $f_{\tau,a,b}$.  First
of all, the following statement is an immediate consequence of
Corollary~\ref{cor:refined} applied to $h(x) = x+a\sin(2\pi x)$ and Fubini's
Theorem:
\begin{cor}
  \label{c.arnoldtongues} Given any $a \in (0,1/2\pi)$, there exists a
  constant $d_0 = d_0(a)$, such that for any $d\geq d_0$ there exists a set
  $\Omega\ssq \kreis$ of positive measure with the following property: 

  For any $\omega \in \Omega$, then there exists a set of positive measure
  $B_\omega \ssq \R$, such that for any $b\in B_\omega$ the qpf Arnold circle
  map $f_{0,a,b}$ is minimal and has a SNA and a SNR.
\end{cor}
 Since we want to study the dependence of the first Arnold tongue on the
parameter $b$, the following notation will be convenient:
\begin{equation}
  A^a_\rho \ := \ \{ (\tau,b) \in \kreis \times \R \mid \rho(f_{\tau,a,b}) =
  \rho \} \ .
\end{equation}
As the rotation number depends monotonically on the parameter $\tau$, there
exist functions $\tau_{a,\rho}^-,\tau_{a,\rho}^+ : \R \ra \kreis$,
such that
\begin{equation}
  A^a_\rho \ = \ \{ (\tau,b) \in \kreis \times \R \mid \tau \in
  [\tau^-_{a,\rho}(b),\tau^+_{a,\rho}(b)] \} \ .
\end{equation}
These functions $\tau^\pm_{a,\rho}$ are continuous for all $a,\rho$ and
coincide (meaning $\tau^-_{a,\rho} =\tau^+_{a,\rho}$) whenever $\rho$ does not
depend rationally on $\omega$, i.e.\ $\rho \notin \Q + \Q\omega \bmod 1$
\cite{bjerkloev/jaeger:2006a}.
\medskip

The canonical lift of the qpf Arnold circle map is given by $$F_{\tau,a,b} \ :
\ \kreis \times \R \xhookleftarrow{} \quad , \quad \thx \mapsto
(\thom,x+\tau+a\sin(2\pi x) + b \cos(2\pi\theta)^d) \ .$$ Obviously there
holds $F_{0,a,b,\theta}(-x) = -F_{0,a,b,\theta}(x)$. This symmetry immediately
implies $\rho(F_{0,a,b}) = 0$, and therefore $0 \in
[\tau^-_{a,0}(b),\tau^+_{a,0}(b)] \ \forall b \in \R$.  On the other hand, if
$d\geq d_0(a)$ and $b\in B_\omega$, where $d_0(a)$ and $B_\omega$ are chosen
as in the above corollary, then $f_{0,a,b}$ is minimal. It therefore follows
from Proposition~\ref{p.strict-monotonicity} that the mapping $\tau \mapsto
\rho(f_{\tau,a,b})$ is strictly monotone at $\tau = 0$.  Consequently, the
first Arnold tongue is collapsed to a single point at this $b$-value, meaning
$\tau^-_{a,0}(b) = \tau^+_{a,0}(b) = 0$. As this happens on a set of $b$ of
positive measure, and since the first Arnold tongue is clearly not collapsed at
$b=0$, the dependence of $\tau^{\pm}_{a,0}$ on $b$ cannot be real-analytic. We
summarise our observations in the following
\begin{prop}
  Suppose $a \in (0,1/2\pi)$ is fixed and $d\geq d_0(a)$, where $d_0(a)$ is
  the constant provided by Corollary~\ref{c.arnoldtongues}~. Let $B_\omega$ be
  as in the corollary. 

Then for any $b\in B_\omega$, there holds $\tau^-_{a,0}(b) =
  \tau^+_{a,0}(b)=0$. Furthermore, the mappings $b \mapsto \tau^\pm_{a,0}(b)$
  are not real-analytic.
\end{prop}
Of course, this raises the question whether the dependence of the boundaries
of the Arnold tongues is analytic in $a$. We have to leave this open here.
However, by the same arguments applied with the roles of $a$ and $b$
interchanged, one obtains the existence of parameters $b$, such that for a set
of $a$'s of positive measure the first Arnold tongue is collapsed. Hence, if
such a parameter $b$ is fixed and the dependence on $a$ was real-analytic,
then the first tongue would have to be reduced to a single point for all $a
\in [0,1/2\pi]$.


\section{Creation of SNA: The basic mechanism} \label{BasicMechanism}

The aim of this section is to prove Theorem~\ref{thm:basic}~. Thereby, we
proceed in three steps. First, we place certain {\em `imaginary'} conditions
of the rotation number $\omega$, and show that these imply the existence of a
sink-source-orbit (Sections~\ref{CriticalSetsDef}
and~\ref{Sink-Source-Orbits}). After this, it remains to show that there exist
rotation numbers which satisfy these conditions. In order to do so, we first
describe the geometry of certain critical sets, which were used before in the
formulation of the conditions on $\omega$ (Section~\ref{CriticalSets}). Using
the obtained information, we then perform a parameter exclusion, which still
leaves a set of positive measure of {\em `good'} $\omega$'s, which have all
the required properties. The technical statements for the parameter exclusion
are contained in Section~\ref{GoodFrequencies}, the final step in the proof is
then given in Section~\ref{BasicProof}~. The proof of the minimality statement
is contained in Section~\ref{Minimality}~.

\subsection{Critical sets and good frequencies} \label{CriticalSetsDef}

\emph{Critical sets.} First we have to define a sequence of \emph{critical
sets}, which project down to \emph{critical regions} and play a major part in
all that follows:
\begin{definition} \label{def:In}
For $\omega \in \kreis$, ${\cal I}_0$ as above and any monotonically
increasing sequence \nofolge{M_n} of integers with $M_0 \geq 2$ we inductively
define nested sequences ${\cal C}_0 ,{\cal C}_1, \ldots$ of \emph{critical
sets} and ${\cal I}_0 \supseteq {\cal I}_1 \supseteq {\cal I}_2 \ldots$ of
\emph{critical regions} in the following way: If ${\cal I}_0 \ld {\cal I}_n$
have been defined, let
\begin{eqnarray*} 
  {\cal A}_n & := & \{\thx \mid \theta \in {\cal I}_n-(M_n-1)\omega,\ x \in C\} \ , \\
  {\cal B}_n & := & \{\thx \mid \theta \in {\cal I}_n+(M_n+1)\omega,\ x \in
  E\} \ , \\
  {\cal C}_n &:= & f^{M_n-1}({\cal A}_n) \cap f^{-M_n-1}({\cal B}_n)
\end{eqnarray*} 
and
\begin{equation} \nonumber
  \label{eq:defIn}
{\cal I}_{n+1} \ := \ \mathrm{int}(\pi_1({\cal C}_n)) \ .
\end{equation}
\end{definition}

\ \\ \emph{Good frequencies.} Further, we impose certain `Diophantine'
conditions on the frequency $\omega$, which mainly state that the critical
sets do not return too fast:
\begin{definition} \label{def:Fn}
Suppose \nofolge{M_n} and \nofolge{{\cal I}_n} are chosen as above and let
\nofolge{K_n} be a monotonically increasing sequence of positive
integers. Further, let \nofolge{\eps_n} be a non-increasing sequence of
positive real numbers which satisfy $\eps_n \geq 3\eps_{n+1} \ \forall
n\in\N_0$. Finally, let
\begin{equation} \nonumber
  {\cal X}_n \ := \ \bigcup_{k=1}^{2K_n M_n} ({\cal I}_n+k\omega) \ \
  \textrm{ and } \ \ {\cal Y}_n \ := \ \bigcup_{j=0}^n \bigcup_{k = -M_j+1}^{M_j+1}
  ({\cal I}_j + k\omega) \ .
\end{equation}
Then we define ${\cal F}_n = {\cal F}_n(M_0 \ld M_n)$ as the
set of those frequencies $\omega \in \kreis$ which satisfy
\begin{equation} \tag*{$({\cal F}1)_n$} \label{eq:F1n}
    d({\cal I}_j,{\cal X}_j) \ > \
    3\eps_j \ \ \ \ \ \forall j = 0 \ld n
\end{equation}
and 
\begin{equation} \tag*{$({\cal F}2)_n$} \label{eq:F2n}
    d(({\cal I}_j -(M_j-1)\omega) \cup ({\cal I}_j+(M_j+1)\omega),{\cal Y}_{j-1}) \ >
    \ 0 \ \ \ \ \ \forall j=1\ld n \ .
\end{equation}
Further, let
\begin{equation} \nonumber
 {\cal Z}_n \ := \ \bigcup_{j=0}^n \bigcup_{k = -M_j+2}^{M_j} ({\cal
 I}_j+k\omega) \ ,
\end{equation}
${\cal Z}_{-1}:= \emptyset$ and ${\cal F}_{-1} := \kreis$.
\end{definition}
Finally let
\[
{\cal V}_{n} \ := \ \bigcup_{j=0}^n \bigcup_{k=1}^{M_j+1} ({\cal I}_j +
k\omega) \ \ \ \textrm{and} \ \ \ {\cal W}_n := \bigcup_{j=0}^n
\bigcup_{k=-M_j+1}^{0} ({\cal I}_j+k\omega) \ 
\]
and ${\cal V}_{-1} = {\cal W}_{-1} = \emptyset$. 
\begin{bem}
  For an easier reading of the following sections, the reader should keep in
  mind the following `intuitive' description of the relations between the sets
  ${\cal Y}_n$, ${\cal Z}_n$, ${\cal V}_n$ and ${\cal W}_n$: ${\cal V}_n$ and
  ${\cal W}_n$ are just the `right' and `left' part of ${\cal Y}_n$, whereas
  ${\cal Z}_n$ is just reduced by one iterate on either side in comparison
  with ${\cal Y}_n$, such that ${\cal Z}_n \pm \omega$ is still contained in
  ${\cal Y}_n$.
\end{bem}


\subsection{Construction of the sink-source-orbits} \label{Sink-Source-Orbits}
 
Recall that for any given point $(\theta_0,x_0)$, we denote its orbit by
$(\theta_k,x_k) = f^k(\theta_0,x_0)$.
\begin{lem}   \label{lem:keyestimateI}
Suppose (\ref{eq:Cinvariance}) holds. Then for all $n \geq 0$, the following
are true:

\ \\
\underline{Forwards iteration:}  \ If
\begin{equation} \label{eq:B1n}
  \tag*{$({\cal B}1)_n$} \left\{ \begin{array}{rcl} \omega & \in & {\cal F}_{n-1}
  \\ \theta_0 & \notin & {\cal Z}_{n-1} \\ x_0 & \in & C \end{array} \right.
\end{equation}
and ${\cal L} \geq 0$ is the first integer, such that $\theta_{{\cal L}} \in {\cal I}_n$, then 
\begin{equation} \label{eq:C1n}
  \tag*{$({\cal C}1)_n$} x_m \notin C \ \follows \ \theta_m \in
  {\cal V}_{n-1} \ \ \ \ \ \forall m=1 \ld {\cal L} \ .
\end{equation}

\ \\
\underline{Backwards iteration:} \ If
\begin{equation} \label{eq:B2n}
  \tag*{$({\cal B}2)_n$} \left\{ \begin{array}{rcl} \omega & \in & {\cal
  F}_{n-1} \\ \theta_0 & \notin & {\cal Z}_{n-1} \\ x_0 & \in & E
  \end{array} \right.
\end{equation}
and ${\cal R}\geq 0$ is the first integer, such that $\theta_{-{\cal R}} \in {\cal I}_n+\omega$,
then 
\begin{equation} \label{eq:C2n}
  \tag*{$({\cal C}2)_n$} x_{-m} \notin E \ \follows \
  \theta_{-m} \in {\cal W}_{n-1} \ \ \
  \ \ \forall m=1 \ld {\cal R} \ .
\end{equation}
\end{lem}
\proof  First of all, note that $({\cal
C}1)_0$ follows directly from (\ref{eq:Cinvariance}).  Now suppose that
\ref{eq:B1n} implies \ref{eq:C1n} and fix $\omega \in {\cal F}_n$, $\theta_0
\notin {\cal Z}_n$ and $x_0 \in C$. Using \ref{eq:F1n} and \ref{eq:F2n},
it is easy to see that
\begin{eqnarray}
  ({\cal I}_n - (M_n-1)\omega) \cap {\cal V}_n & = & \emptyset \ , \label{eq:In-Ln:Vn}\\
  ({\cal I}_n + (M_n+1)\omega) \cap {\cal I}_0 & = & \emptyset \ ,\label{eq:In+Rn:I0} \\
  ({\cal I}_n + (M_n+2)\omega) \cap {\cal Z}_n & = & \emptyset  \ \label{eq:In+Rn+2:Zn} .
\end{eqnarray}
Let ${\cal L}$ be the first integer such that $\theta_{\cal L} \in {\cal
  I}_{n+1}$ and let $0 < L_1 < L_2 < \ldots < L_J = {\cal L}$ be those times
$0\leq i \leq {\cal L}$ with $\theta_i \in {\cal I}_n$. If we denote condition
$({\cal C}1)_{n+1}$ with ${\cal L}$ replaced by $L_j$ by $({\cal
  C}1)_{n+1}[L_j]$, then $({\cal C}1)_{n+1}[L_1]$ follows from \ref{eq:C1n}
(note that ${\cal Z}_{n-1} \ssq {\cal Z}_n$, ${\cal F}_n \ssq {\cal F}_{n-1}$
and ${\cal V}_{n-1} \ssq {\cal V}_n$).

Assume now that $({\cal C}1)_{n+1}[L_j]$ holds for some $1 \leq j < J$. As
$\theta_0 \notin {\cal Z}_n$ we have $L_j - M_n +1\geq 0$, and as
$\theta_{L_j-M_n+1} \notin {\cal V}_n$ due to (\ref{eq:In-Ln:Vn}) it follows that
$x_{L_j-M_n+1} \in C$. Consequently
$(\theta_{L_j-M_n+1},x_{L_j-M_n+1}) \in {\cal A}_n$, and as $\theta_{L_j} \notin
{\cal I}_{n+1}$ we must have $(\theta_{L_j+M_n+1},x_{L_j+M_n+1}) \notin {\cal B}_n$,
which means
\[
    x_{L_j+M_n+1} \ \notin \ E \ .
\]
As $\theta_{L_j+M_n+1} \notin {\cal I}_0$ by (\ref{eq:In+Rn:I0}) we can
apply (\ref{eq:Cinvariance}) and obtain $x_{L_j+M_n+2} \in
C$. Before, we could have had $x_k \notin C$ for some $k
\in \{L_j+1 \ld L_j+M_n+1\}$, but for such $k$ there obviously holds
\[
  \theta_k \ \in \ {\cal I}_n+\omega \cup \ldots \cup {\cal I}_n+(M_n+1)\omega \
  \ssq \ {\cal V}_n \ .
\] 
Further, as $\theta_{L_j+M_n+2} \notin {\cal Z}_n\supseteq {\cal Z}_{n-1}$ by
(\ref{eq:In+Rn+2:Zn}) and ${\cal F}_n \ssq {\cal F}_{n-1}$, we can now apply
\ref{eq:C1n} and obtain $({\cal C}1)_{n+1}[L_{j+1}]$. As $L_J = {\cal L}$, this
completes the proof of $({\cal C}1)_{n+1}$. 

{ \ \\ \underline{\emph{Backwards iteration:}} $({\cal C}2)_0$
follows directly from (\ref{eq:Einvariance}). Suppose that \ref{eq:B1n}
implies \ref{eq:C2n} and fix $\omega \in {\cal F}_n$, $\theta_0 \notin {\cal Z}_n$
and $x_0 \in E$. Using \ref{eq:F1n} and \ref{eq:F2n}, we see that
\begin{eqnarray}
  ({\cal I}_n + (M_n+1)\omega) \cap {\cal W}_n & = & \emptyset \ ,
  \label{eq:In+Rn:Wn}\\ ({\cal I}_n - (M_n-1)\omega) \cap {\cal I}_0 & = &
  \emptyset \ ,\label{eq:In-Ln:I0} \\ ({\cal I}_n - M_n\omega) \cap {\cal
  Z}_n & = &
  \emptyset \ \label{eq:In-Ln:Zn} .
\end{eqnarray}
Let ${\cal R}$ be the first integer such that $\theta_{-{\cal R}} \in {\cal
  I}_{n+1}+\omega$ and let $0 < R_1 < R_2 < \ldots < R_J = {\cal R}$ be those times $0\leq i \leq {\cal R}$
with $\theta_{-i} \in {\cal I}_n+\omega$. If we denote condition $({\cal
C}2)_{n+1}$ with ${\cal R}$ replaced by $R_j$ by $({\cal C}2)_{n+1}[R_j]$, then
$({\cal C}2)_{n+1}[R_1]$ follows from \ref{eq:C2n}.

Assume now that $({\cal C}2)_{n+1}[R_j]$ holds for some $1 \leq j < J$. As
$\theta_0 \notin {\cal Z}_n$ we have $R_j - M_n \geq 0$, and as
$\theta_{-R_j+M_n} \notin {\cal W}_n$ due to (\ref{eq:In+Rn:Wn}) (note that
$\theta_{-R_j-1} \in {\cal I}_n$) it follows that $x_{-R_j+M_n} \in
E$. Consequently $(\theta_{-R_j+M_n},x_{-R_j+M_n}) \in {\cal B}_n$, and as
$\theta_{-R_j-1} \notin {\cal I}_{n+1}$ we must have
$(\theta_{-R_j-M_n},$ $x_{-R_j-M_n})$ $\notin {\cal A}_n$, which means
\[
    x_{-R_j-M_n} \ \notin \ C \ .
\]
As $\theta_{-R_j-M_n} \notin {\cal I}_0$ by (\ref{eq:In-Ln:I0}) we can
apply (\ref{eq:Einvariance}) and obtain $x_{-R_j-M_n-1} \in
E$. Before, we could have had $x_{-k} \notin E$ for some
$k \in \{R_j+1 \ld R_j+M_n\}$, but for such $k$ there obviously holds
\[
  \theta_k \ \in \ {\cal I}_n \cup {\cal I}_n - \omega \cup \ldots \cup {\cal
  I}_n-M_n\omega \ \ssq \ {\cal W}_n \ .
\] 
Further, as $\theta_{-R_j-M_n-1} \notin {\cal Z}_n\supseteq {\cal Z}_{n-1}$ by
(\ref{eq:In-Ln:Zn}) and ${\cal F}_n \ssq {\cal F}_{n-1}$, we can now apply
\ref{eq:C2n} and obtain $({\cal C}2)_{n+1}[R_{j+1}]$.} As $R_J = {\cal R}$,
this completes the proof.

\qed

\begin{bem} \label{bem:keyestimate}
 \alphlist 
 \item Suppose (\ref{eq:Cinvariance}) holds, $\omega \in {\cal F}_n$ and
  $(\theta_0,x_0) \in {\cal A}_n$. Then \ref{eq:B1n} holds and ${\cal
  L}=M_n-1$.

 In order to see this, note that $x_0 \in C$ holds by definition
  of ${\cal A}_n$, and $\theta_0 \notin {\cal Z}_{n-1}$ follows from
  \begin{eqnarray} \label{eq:In-Ln+1:Zn}
  ({\cal I}_n - (M_n-1)\omega) \cap  {\cal Z}_n & = & \emptyset \ ,
  \end{eqnarray}
  which is a consequence of see \ref{eq:F1n} and \ref{eq:F2n}.
\item 
  Similarly, suppose $\omega\in {\cal F}_{n-1}$ and $(\theta_0,x_0) \in {\cal
  B}_n$. Then \ref{eq:B2n} holds and ${\cal R}=M_n$. 

This follows by the same argument as (a): $x_0 \in E$ holds by definition of ${\cal
  B}_n$ and $\theta_0\notin {\cal Z}_{n-1}$ follows from
  \begin{equation}
  \label{eq:In+Rn:Zn} ({\cal I}_n + (M_n+1)\omega) \cap {\cal Z}_n \ = \
  \emptyset \ ,
  \end{equation}
  which is again a consequence of \ref{eq:F1n} and \ref{eq:F2n}.  \listend
\end{bem}

\begin{cor}
  \label{cor:keyestimate}
Suppose  (\ref{eq:Cinvariance}) holds and $\omega \in {\cal F}_n$. Then
\begin{equation} \label{eq:AnBninclusions}
  f^{M_n-M_{n-1}}({\cal A}_n) \ \ssq \ {\cal A}_{n-1} \ \ \ \textrm{and} \ \ \ 
  f^{-M_n+M_{n-1}}({\cal B}_n) \ \ssq \ {\cal B}_{n-1} \ .
\end{equation}
Consequently ${\cal C}_0 \supseteq {\cal C}_1 \supseteq {\cal C}_2 \supseteq \ldots$~. Further
\begin{equation} \label{eq:laststepinclusion}
  f^{M_n-1}({\cal A}_n) \ \ssq \ {\cal I}_n \times C \ \ \ \textrm{and} \ \
  \ f^{-M_n}({\cal B}_n) \ \ssq \ ({\cal I}_n+\omega) \times E \ .
\end{equation} 
\end{cor}
\proof Let $(\theta_0,x_0) \in {\cal A}_n$, such that, by the preceding
remark, \ref{eq:B1n} holds. There holds $$({\cal I}_n - (M_{n-1}-1)\omega)
\cap {\cal V}_{n-1} = \emptyset \ .$$ This follows from (\ref{eq:In-Ln:Vn}),
applied to $n-1$ and using that ${\cal I}_n \ssq {\cal I}_{n-1}$. Therefore we
have $\theta_{M_n-M_{n-1}} \notin {\cal V}_{n-1}$, such that we can apply
Lemma~\ref{lem:keyestimateI} and obtain that $x_{M_n-M_{n-1}} \in C$, which
means that $f^{M_n-M_{n-1}}(\theta_0,x_0) \in {\cal A}_{n-1}$. As
$(\theta_0,x_0) \in {\cal A}_n$ was arbitrary, this proves the first inclusion
in (\ref{eq:AnBninclusions}), and the argument for the second one is similar.
Finally, as
\begin{equation} \label{eq:In:Vn}
  {\cal I}_n \cap {\cal V}_n \ = \ \emptyset \ 
\end{equation}
and
\begin{equation} \label{eq:In+1:Wn}
  {\cal I}_n +\omega \cap {\cal W}_n \ = \ \emptyset \
\end{equation}
due to \ref{eq:F1n}, the inclusions in (\ref{eq:laststepinclusion}) follow in
the same way.

\qed
\bigskip

The preceding lemma gives some first control about the time an orbit
spends in the expanding and contracting region. In order to make use of this
information, we need to quantify it. For given $\omega,\theta_0, x_0$ and
$0\leq m \leq N$ let
\begin{eqnarray} \label{eq:PNm}
{\cal P}^N_m & := & \#\{k\in[m,N-1] \mid x_k \in C \} \ , \\ \label{eq:QNm}
{\cal Q}^N_m & := & \#\{k\in[m,N-1] \mid x_{-k} \in E\} \ .
\end{eqnarray}
Further, let $\beta_0 = 1$ and 
\begin{equation}
\beta_n \ : =  \ \prod_{j=0}^{n-1} \left( 1-\frac{1}{K_j} \right) \ .
\end{equation}
\begin{lem} \label{lem:occupationtimes}
Suppose  (\ref{eq:Cinvariance}) holds. Then for all $n\geq 0$ the following
are true:

\ \\ \underline{Forwards iteration:} \ Suppose \ref{eq:B1n} holds and let
  ${\cal L}$ be chosen as in Lemma~\ref{lem:keyestimateI}~.  Then
\begin{equation} \label{eq:C3n} \tag*{$({\cal C}3)_n$}
  {\cal P}_m^{{\cal L}} \ \geq \ \beta_n \cdot ({\cal L}-m) \ \ \ \ \
  \forall m = 0 \ld {\cal L}-1 \ .
\end{equation}
Further $x_{\cal L} \in C$. 

\ \\ \underline{Backwards iteration:} \ Suppose \ref{eq:B2n} holds and let
${\cal R}$ be chosen as in Lemma~\ref{lem:keyestimateI}~.  Then
\begin{equation} \label{eq:C4n} \tag*{$({\cal C}4)_n$} {\cal Q}_m^{{\cal R}} \
  \geq \ \beta_n \cdot ({\cal R}-m) \ \ \ \ \ \forall m = 0 \ld {\cal R}-1 \ .
\end{equation}
Further $x_{-{\cal R}} \in E$.
\end{lem}
\proof As ${\cal V}_{-1}$ is void, $({\cal C}3)_0$ follows directly from
$({\cal C}1)_0$. Suppose that \ref{eq:B1n} implies \ref{eq:C3n} and fix
$\omega \in {\cal F}_n$, $\theta_0 \notin {\cal Z}_n$ and $x_0 \in C$. As in
the proof of Lemma~\ref{lem:keyestimateI}, let $0 < L_1 < L_2 < \ldots < L_J =
{\cal L}$ be those times $0\leq i \leq {\cal L}$ with $\theta_i \in {\cal
I}_n$ and denote condition \ref{eq:C3n} with ${\cal L}$ replaced by $L_j$ by
$({\cal C}3)_n[L_j]$.

As $\beta_{n+1} \leq \beta_n$, condition $({\cal C}3)_{n+1}[L_1]$ follows
  from \ref{eq:C3n}. Suppose $({\cal C}3)_{n+1}[L_j]$ holds for some $1 \leq j
  <J$. Using \ref{eq:F1n} and \ref{eq:F2n} we see that
  \begin{equation}
    \label{eq:In+Rn+2:Vn}
   ({\cal I}_n + (M_n+2)\omega) \cap {\cal V}_n \ = \ \emptyset \ ,
  \end{equation}
such that in particular $\theta_{L_j+M_n+2} \notin {\cal V}_n$ and
consequently $x_{L_j+M_n+2} \in C$ by $({\cal C}1)_{n+1}$. As further
$\theta_{L_j+M_n+2} \notin {\cal Z}_n$ by (\ref{eq:In+Rn+2:Zn}), \ref{eq:C3n}
implies that for any $ m \in [L_j+M_n+2,L_{j+1}]$ there holds
\begin{equation} \label{eq:Pm1}
    {\cal P}_m^{L_{j+1}} \ \geq \ \beta_n \cdot (L_{j+1}-m) \ .
\end{equation}
This proves $({\cal C}3)_{n+1}[L_{j+1}]$ for such $m$. Further, by \ref{eq:F1n}
  we have $L_{j+1} - L_j > 2 K_n M_n$. Hence, for any $m \in
  [L_j,L_j+M_n+1]$ we obtain the estimate
\begin{eqnarray*}
{\cal P}_m^{L_{j+1}} & \geq & {\cal P}_{L_j+M_n+2}^{L_{j+1}} \ \geq \
\beta_n \cdot (L_{j+1} - L_j-M_n-2) \\ & \geq & \beta_n \cdot
\frac{L_{j+1} - L_j - M_n-2}{L_{j+1} - L_j} \cdot (L_{j+1}-m) \\ & \geq &
\beta_n \cdot \left(1-\frac{M_n+2}{2K_nM_n} \right) (L_{j+1}-m) \
\stackrel{M_0 \geq 2}{\geq} \ \beta_{n+1} \cdot
(L_{j+1}-m) \ .
\end{eqnarray*}
Finally, if $m\in [0,L_j]$ the statement follows by combining the estimate for
${\cal P}^{L_{j+1}}_{L_j}$ with the one for ${\cal P}^{L_j}_m$ obtained from $({\cal
C}3)_{n+1}[L_j]$.

{\ \\ \underline{\emph{Backwards iteration:}} \ As ${\cal W}_{-1}$ is void,
$({\cal C}4)_0$ follows directly from $({\cal C}2)_0$.  Suppose that
\ref{eq:B2n} implies \ref{eq:C4n} and fix $\omega \in {\cal F}_n$, $\theta_0
\notin {\cal Z}_n$ and $x_0 \in E$. Let $0 < R_1 < R_2 < \ldots < R_J = {\cal
R}$ be those times $0\leq i \leq {\cal L}$ with $\theta_{-i} \in {\cal
I}_n+\omega$ and denote condition \ref{eq:C4n} with ${\cal R}$ replaced by
$R_j$ by $({\cal C}4)_n[R_j]$.

As $\beta_{n+1} \leq \beta_n$, condition $({\cal C}4)_{n+1}[R_1]$ follows
  from \ref{eq:C4n}. Suppose $({\cal C}4)_{n+1}[R_j]$ holds for some $1 \leq j
  <J$. Using \ref{eq:F1n} and \ref{eq:F2n} we see that
  \begin{equation}
    \label{eq:In-Ln:Wn}
   ({\cal I}_n - M_n\omega) \cap {\cal W}_n \ = \ \emptyset \ ,
  \end{equation}
such that in particular $\theta_{-R_j-M_n-1} \notin {\cal W}_n$ and
consequently $x_{-R_j-M_n-1} \in E$ by $({\cal C}2)_{n+1}$. As
further $\theta_{-R_j-M_n-1} \notin {\cal Z}_n$ by (\ref{eq:In-Ln:Zn}),
\ref{eq:C4n} implies that for any $ m \in [R_j+M_n+1,R_{j+1}]$ there holds
\begin{equation} \label{eq:Qm1}
    {\cal Q}_m^{R_{j+1}} \ \geq \ \beta_n \cdot (R_{j+1}-m) \ .
\end{equation}
This proves $({\cal C}4)_{n+1}[R_{j+1}]$ for such $m$. Further, by \ref{eq:F1n}
  we have $R_{j+1} - R_j  >2 K_n M_n$. Hence, for any $m \in
  [R_j,R_j+L_n]$ we obtain the estimate
\begin{eqnarray*}
{\cal Q}_m^{R_{j+1}} & \geq & {\cal Q}_{R_j+M_n+1}^{R_{j+1}} \ \geq \ \beta_n
\cdot (R_{j+1} - R_j -M_n-1) \\ & \geq & \beta_n \cdot \frac{R_{j+1} - R_j -
M_n -1}{R_{j+1} - R_j} \cdot (R_{j+1}-m) \\ & \geq & \beta_n \cdot
\left(1-\frac{M_n+1}{2K_nM_n} \right) (R_{j+1}-m) \ \geq
\ \beta_{n+1} \cdot (R_{j+1}-m) \ .
\end{eqnarray*}
Finally, if $m\in [0,R_j]$ the statement follows by combining the estimate
for ${\cal P}^{R_{j+1}}_{R_j}$ with the one that $({\cal C}4)_{n+1}[R_j]$
yields for ${\cal P}^{R_j}_m$.}

\qed

\ \\ 
Let
\begin{equation}
\beta \ := \  \nLim \beta_n \ = \ \inf_n \beta_n \
\end{equation}
and
\begin{equation}
\alpha_-  :=  \alpha_c^{\beta}  \alpha_u^{1-\beta} \  ,\quad
\alpha_+  :=  \alpha_e^{\beta}  \alpha_l^{1-\beta} \ .
\end{equation}

\begin{cor}
  \label{cor:derivatives} Suppose  (\ref{eq:Cinvariance})--(\ref{eq:bounds3})
  hold and $\omega \in {\cal F}_n$. If $(\theta,x)
  \in \mbox{cl}(f^{M_n}({\cal A}_n))$, then for all $k \in [0,M_n]$ there holds
  \begin{equation} \label{eq:derivativesI}
    \partial_x f_\theta^{-k}(x) \ \geq \ \alpha_-^{-k} \ .
  \end{equation}
If $(\theta,x) \in \mbox{cl}(f^{-M_n}({\cal B}_n))$, then for all $k
\in [0,M_n]$ there holds
\begin{equation} \label{eq:derivativesII}
  \partial_x f^k_\theta(x) \ \geq \  \alpha_+^k \ .
\end{equation}
\end{cor}
\proof By continuity, it suffices to prove the above estimates on
$f^{M_n}({\cal A}_n)$ and $f^{-M_n}({\cal B}_n)$, respectively. We start by
proving (\ref{eq:derivativesII}). 

Suppose $(\theta,x) \in f^{-M_n}({\cal B}_n)$ and let $(\theta_0,x_0) =
f^{M_n}(\theta,x) \in {\cal B}_n$. Then due to Remark~\ref{bem:keyestimate} we
have that ${\cal R} = M_n$ and \ref{eq:B2n} holds. Using (\ref{eq:bounds1}),
(\ref{eq:bounds2}) and the fact that $x = x_{-{\cal R}} \in E$ (see
Lemma~\ref{lem:occupationtimes}) we obtain
\begin{eqnarray} \label{eq:factorcount}
  \partial_x f^k_\theta(x) & = & \prod_{j={\cal R}-k+1}^{\cal R} \partial_x
  f_{\theta_{-j}}(x_{-j}) \ \geq \ \alpha_e \cdot \alpha_e^{{\cal Q}^{\cal
  R}_{{\cal R}-k+1}} \cdot \alpha_l^{k-1-{\cal Q}^{\cal R}_{{\cal R}-k+1}}
\end{eqnarray}
Applying \ref{eq:C4n} and using that $\alpha_e \geq \alpha_+$ yields the statement.

As $\partial_x f_\theta^{-k}(x) = (\partial_x
f_{\theta-k\omega}(f_\theta^{-k}(x)))^{-1}$, the estimate in
(\ref{eq:derivativesI}) can be obtained in the same way.

\qed

\begin{prop}
  \label{prop:firststep}
Suppose  (\ref{eq:Cinvariance})--(\ref{eq:bounds3})
  hold, $\min\{\alpha_-^{-1},\alpha_+\} > 1$, $\omega \in \ncap {\cal F}_n$ and all critical sets
${\cal I}_n$ are non-void. Then $f$ has a sink-source-orbit.
\end{prop}
\proof As all critical sets ${\cal I}_n$ are non-void, the same is obviously
true for the sets $\mbox{cl}({\cal C}_n)=$ and their images $\mbox{cl}(f({\cal
C}_n))=\mbox{cl}(f^{M_n}({\cal A}_n))\cap\mbox{cl}(f^{-M_n}({\cal B}_n))$. Due
to Corollary~\ref{cor:keyestimate}, the later form a nested sequence of
compact sets, such that their intersection is non-void as well. Let
$(\theta,x) \in \ncap \mbox{cl}(f({\cal C}_n))$. Then due to
(\ref{eq:derivativesI}) and as $M_n \nearrow \infty$, we obtain
\[
    \lambda^-(\theta,x) = \limsup_{k\ra \infty} \ktel \log |\partial_x
    f^{-k}_\theta(x)| \geq -\log \alpha_- > 0 \ .
\]
and similarly (\ref{eq:derivativesII}) yields $\lambda^+(\theta,x) \geq \log
\alpha_+ > 0$. 

\qed


\subsection{Geometry of the critical sets} \label{CriticalSets}

In this section we turn to the description of the critical sets ${\cal C}_n$
and the corresponding critical regions ${\cal I}_{n+1}$. In particular, we
want to obtain information about their size and their dependence on $\omega$
(which we have kept implicit so far). Suppose
$I=I(\omega)=(a(\omega),b(\omega))$ is a connected component of ${\cal
I}_n$. Then we use the notation
\[
  |\partial_\omega I| \ = \ \max\{|\partial_\omega a(\omega)|,|\partial_\omega
  b(\omega)|\} \ ,
\]
provided both derivatives on the right side exist. In this case we call $I$
differentiable with respect to $\omega$. We will use the following
inductive assumption:
\begin{equation} \label{eq:In} \tag*{$({\cal I})_n$}
  \left\{ \begin{array}{cl} (i) & \textrm{For each } j\in[0,n] , \ {\cal I}_j
    \textrm{ consists of } {\cal N} \textrm{ disjoint open intervals } \\&
    I^1_j \ld I^{\cal N}_j. \\ \\ (ii) & \textrm{For } j \in [1,n], \textrm{
    each connected component of } {\cal I}_{j-1} \textrm{ contains} \\ &
    \textrm{exactly one connected component of } {\cal I}_j. \textrm{ Thus, by
    suitable } \\ & \textrm{labelling, } I^\iota_j \ssq I^\iota_{j-1} \
    \forall \iota = 1 \ld {\cal N}. \\ \\ (iii)& \textrm{For all } j \in [0,n]
    \textrm{ the set } {\cal F}_j \textrm{ is open and all } I^\iota_j
    \textrm{ are differentiable} \\ & \textrm{ with respect to } \omega
    \textrm{ on } {\cal F}_j.
  \end{array} \right. 
\end{equation}
Note $({\cal I})_0$ follows directly from the choice of ${\cal I}_0$ in
Section~\ref{MainAssumptions} and the definition of ${\cal F}_0$. (The second
statement is void for $n=0$.)
\bigskip

In order to describe the geometry of the critical sets ${\cal C}_n$, or rather
their images $f({\cal C}_n)$, we have to introduce some notation and make some
preliminary remarks, which we will use in the whole section. For any $\iota \in
[1,{\cal N}]$ we let
  \begin{eqnarray}
  {\cal A}_n^\iota & := & \{\thx \mid \theta \in I^\iota_n-(M_n-1)\omega,\ x \in C\} \ ,
  \\ {\cal B}_n^\iota & := & \{\thx \mid \theta \in I_n^\iota+(M_n+1)\omega,\ x \in E\} \
  .
  \end{eqnarray}
For $\theta \in {\cal I}_n+\omega$ let 
\begin{equation} \label{eq:phipsidef}
\varphi_{\iota,n}^\pm(\theta) \ := \ f_{\theta-M_n\omega}^{M_n}(c^\pm)
\quad \textrm{ and } \quad \psi_{\iota,n}^\pm(\theta) \ := \
f_{\theta+M_n\omega}^{-M_n}(e^\pm)\ ,
\end{equation}
such that
\begin{eqnarray*}
  f^{M_n}({\cal A}_n^\iota) & = & \{ \thx \mid \theta \in I_n^\iota+\omega,\
  x \in [\varphi_{\iota,n}^-(\theta) ,\varphi_{\iota,n}^+(\theta)] \} \ , \\
  f^{-M_n}({\cal B}_n^\iota) & = & \{ \thx \mid \theta \in I_n^\iota+\omega,\
  x \in [\psi_{\iota,n}^-(\theta),  \psi_{\iota,n}^+(\theta)] \} \ .
\end{eqnarray*}
In order to start the induction, it is also convenient to define
\begin{equation} \label{eq:phipsi-1-def}
  \varphi^{\pm}_{-1}(\theta) \ := \ f_{\theta-\omega}(c^\pm) \quad \textrm{
  and } \quad \psi^{\pm}_{-1}(\theta) \ := \
  e^\pm \ .
\end{equation}

In all of the proofs of this section we will always fix $\iota$ in order to
concentrate on one connected component of ${\cal I}_n$. In principle we would
have to distinguish two cases, namely that of an upwards and that of a
downwards crossing (see (\ref{eq:crossing})). However, as the two cases are
completely symmetric we can always assume, without loss of generality, that
the crossing between $f^{M_n}({\cal A}_n)$ and $f^{-M_n}({\cal B}_n)$ is
`upwards', that is $\partial_\theta f_\theta(x) > s$ on $I^\iota_n
\subseteq I^\iota_0$.

Then the second inductive assumption which will be used in this section is the
following: Suppose that $I^\iota_n(\omega) =
(a_{\iota,n}(\omega),b_{\iota,n}(\omega))$ and let $J^\varphi_n(\theta) :=
(\varphi^-_n(\theta),\varphi^+_n(\theta))$ and $J_n^\psi(\theta) :=
(\psi^-_n(\theta),\psi_n^+(\theta))$. Then we will assume that
\begin{equation}\tag*{$(\Phi/\Psi)_n$} \label{eq:phipsin}
   \begin{array}{lcrc}
     J^\varphi_{n-1}(a_{\iota,n}(\omega)+\omega) \cap
      J^\psi_{n-1}(a_{\iota,n}(\omega)+\omega) & = & \emptyset \\
     J^\varphi_{n-1}(b_{\iota,n}(\omega)+\omega) \cap
      J^\psi_{n-1}(b_{\iota,n}(\omega)+\omega) & = & \emptyset
    \end{array}
\end{equation}
Note that due to the definition of $\varphi^\pm_{-1}$ and $\psi^{\pm}_{-1}$ in
(\ref{eq:phipsi-1-def}), the statement $(\Phi/\Psi)_0$ is a consequence of
(\ref{eq:Cinvariance}).
\bigskip

 Now we can derive some estimates concerning the geometry of the sets $f({\cal
C}_n)$. We start with an easy one. Let
\begin{eqnarray}
   h^\varphi_n  & := &  \inf_{\theta \in {\cal I}_n+\omega}
  |\varphi_n^+(\theta)-\varphi_n^-(\theta)|\ , \\
  h^\psi_n & := & \inf_{\theta \in {\cal I}_n+\omega}
  |\psi_n^+(\theta)-\psi_n^-(\theta)| \ ,\\
  H^\varphi_n  & := & \sup_{\theta \in {\cal I}_n+\omega}
  |\varphi_n^+(\theta)-\varphi_n^-(\theta)| \ , \\
  H^\psi_n & := & \sup_{\theta \in {\cal I}_n+\omega}
  |\psi_n^+(\theta)-\psi_n^-(\theta)| \ .
\end{eqnarray}
\begin{lem} \label{lem:graphdistance} \ Suppose
  (\ref{eq:Cinvariance})--(\ref{eq:bounds3}) hold and $\omega \in {\cal
  F}_n$. Then
\begin{eqnarray}
\label{eq:phidistance}    |C| \cdot \alpha_l^{M_n} \ \leq &  h^\varphi_n \
  \leq \ H^\varphi_n & \leq \ |C| \cdot \alpha_-^{M_n} 
\end{eqnarray}
and 
\begin{eqnarray} |E| \cdot
  \alpha_u^{-M_n} \ \leq & h^\psi_n \ \leq \ H^\psi_n & \leq \ |E|\cdot
  \alpha_+^{-M_n} .\label{eq:psidistance}
\end{eqnarray}
\end{lem}
\proof As the vertical size of the sets ${\cal A}_n$ and ${\cal B}_n$ is $|C|$
  and $|E|$, respectively, the lower bounds are a direct consequence of
  (\ref{eq:bounds1}) and the upper bounds follow from
  Corollary~\ref{cor:derivatives}~.

\qed
\bigskip

Next, we turn to some more serious estimates. Let
\begin{eqnarray}
  l^\varphi_n & := & \inf_{\theta\in{\cal I}_n+\omega} \left|\partial_\theta
  \varphi_n^\pm(\theta)\right| \ , \\
  u^\varphi_n & := &  \sup_{\theta\in{\cal I}_n+\omega} \left|\partial_\theta
  \varphi_n^\pm(\theta)\right| \ , \\ u^\psi_n & := &
  \sup_{\theta\in{\cal I}_n+\omega} \left|\partial_\theta
  \psi_n^\pm(\theta)\right| \ .
\end{eqnarray}
\begin{lem} \label{lem:dthgraphs}
 Suppose  (\ref{eq:Cinvariance})--(\ref{eq:crossing})
  hold and $\omega \in {\cal F}_n$. Then 
\begin{equation}\label{eq:dthphi}
 s - S/(\alpha_-^{-1}-1)
\ \leq \ l^\varphi_n \ \leq \ u^\varphi_n \ \leq \ S
+ S/(\alpha_-^{-1}-1) 
\end{equation}
and 
\begin{equation}
 \label{eq:dthpsi} 
u^\psi_n \ \leq \ S/(\alpha_+-1) \ .
\end{equation}
\end{lem}
\proof  In order to prove (\ref{eq:dthphi}), note that for any ${\cal L} \in
\N$ and $(\theta_0,x_0) \in \ntorus$ there holds
\begin{equation} \label{eq:dthstepI}
\partial_\theta f^{{\cal L}+1}_{\theta_0}(x_0) \ = \ \partial_\theta
f^{{\cal L}}_{\theta_1}(x_1) +
\partial_xf^{{\cal L}}_{\theta_1}(x_1) \cdot \partial_\theta
f_{\theta_0}(x_0) \ .
\end{equation}
By induction, we thus obtain
\begin{equation} 
\partial_\theta f^{{\cal L}+1}_{\theta_0}(x_0) \ = \ \partial_\theta
f_{\theta_{\cal L}}(x_{\cal L}) \ + \ \sum_{k=0}^{{\cal L}-1} \partial_x
f^{{\cal L}-k}_{\theta_{k+1}}(x_{k+1}) \cdot
\partial_\theta f_{\theta_{k}}(x_{k}) \ .
\label{eq:dthchainI}
\end{equation}
Now suppose $\theta \in I_n^\iota+\omega$ and let $(\theta_0,x_0) =
(\theta-M_n\omega,c^\pm)$ and ${\cal L} = M_n-1$, such that $f^{{\cal
L}+1}_{\theta_0}(x_0) = \varphi_n^\pm(\theta)$. Note that thus ${\cal L}$
coincides with the choice in Lemma~\ref{lem:keyestimateI} (see
Remark~\ref{bem:keyestimate}). By (\ref{eq:bounddth}) and (\ref{eq:s}) we have
\[
s \ < \ \left|\partial_\theta
f_{\theta_{\cal L}}(x_{\cal L})\right| \ < \ S \ .
\]
Further, using (\ref{eq:derivativesI}) from Corollary \ref{cor:derivatives}
we obtain that
\begin{equation} \label{eq:dxestimate}
\left|\partial_x f^{{\cal L}-k}_{\theta_{k+1}}(x_{k+1})\right| \ = \
\left|\left(\partial_x f^{-({\cal L} -k)}_{\theta_{{\cal L}+1}}(x_{{\cal
L}+1}) \right)^{-1}\right| \ \leq \  \alpha_-^{{\cal L}-k} \ .
\end{equation}
As $|\partial_\theta f_{\theta_k}| \leq S \ \forall k$ by
(\ref{eq:bounddth}), this yields the required estimates.

\ \\ The proof of (\ref{eq:dthpsi}) is slightly more intricate. First of all,
similar to (\ref{eq:dthchainI}) we obtain that for any ${\cal R} \in \N$ and
$(\theta_0,x_0) \in \ntorus$
\begin{equation} \label{eq:dthchainII}
\partial_\theta f^{-\cal R}_{\theta_0}(x_0) \ = \ \sum_{k=1}^{\cal R} \partial_x
f^{-{\cal R}+k}_{\theta_{-k}}(x_{-k}) \cdot
\partial_\theta f^{-1}_{\theta_{-k+1}}(x_{-k+1}) \ .
\end{equation}
Let $(\theta_0,x_0) = (\theta+M_n\omega,e^\pm)$ and ${\cal R} = M_n$, such
that $f^{-\cal R}_{\theta_0}(x_0) = \psi_n^\pm(\theta)$. Again, this coincides
with the choice of ${\cal R}$ in Lemma~\ref{lem:keyestimateI}. In order to
obtain an estimate on the second factor in the sum in (\ref{eq:dthchainII}), we
note that 
\[
0 \ = \ \partial_\theta \left(f_{\theta-\omega} \circ f_\theta^{-1}(x)\right) \ = \
\partial_\theta f_{\theta-\omega}(f_{\theta}^{-1}(x)) +
\partial_x f_{\theta-\omega}(f_{\theta}^{-1}(x)) \cdot \partial_\theta
f^{-1}_{\theta}(x) \ ,
\]
such that
 \begin{equation} \label{eq:dthfinbound}
\left|\partial_\theta f^{-1}_{\theta}(x)\right| \ \leq \ \left| \frac{\partial_\theta
 f_{\theta-\omega}(f_{\theta}^{-1}(x))}{\partial_x
 f_{\theta-\omega}(f_\theta^{-1}(x))}\right| \ .
\end{equation}
Therefore $\left|\partial_\theta f^{-1}_{\theta_{-k+1}}(x_{-k+1})\right|$ will
be smaller than $\frac{|\partial_\theta f_{\theta_{-k}}(x_{-k})|}{\alpha_e}$
whenever $x_{-k} \in E$ and always smaller than $\frac{|\partial_\theta
f_{\theta_{-k}}(x_{-k})|}{\alpha_l}$. Combining this with
(\ref{eq:factorcount}) yields
\begin{eqnarray} \nonumber
\lefteqn{\left|\partial_x f^{-{\cal R}+k}_{\theta_{-k}}(x_{-k}) \cdot
\partial_\theta f^{-1}_{\theta_{-k+1}}(x_{-k+1})\right| \ =} \\
\label{eq:dthdxproductbound} & = & \left|\left(\partial_x f^{{\cal
R}-k}_{\theta_{-\cal R}}(x_{-\cal R})\right)^{-1} \cdot \partial_\theta
f^{-1}_{\theta_{-k+1}}(x_{-k+1})\right| \\ \nonumber &\leq & \alpha_e^{-1}
\cdot \alpha_e^{-{\cal Q}^{\cal R}_{k}} \cdot \alpha_l^{-({\cal R}-k-{\cal
Q}^{\cal R}_{k})} \cdot |\partial_\theta f_{\theta_{-k}}(x_{-k})|\\ & \leq &
 \alpha_+^{-({\cal R}+1-k)} \cdot |\partial_\theta
f_{\theta_{-k}}(x_{-k})| \nonumber \ \leq \ 
\alpha_+^{-({\cal R}+1-k)} \cdot S \ ,
\end{eqnarray}
and summing up over $k$ proves~(\ref{eq:dthpsi}). 

\qed
\bigskip

For the remainder of this section, we will write $\varphi^\pm_n(\theta) =
\varphi^\pm_n(\theta,\omega)$ and
$\psi_n^\pm(\theta)=\psi_n^\pm(\theta,\omega)$, in order to make the dependence
on $\omega$ explicit. Let
\begin{eqnarray}
  \gamma_n^\varphi & := & \sup_{\theta \in {\cal I}_n+\omega} \left|\partial_\theta
   \varphi^\pm_n(\theta,\omega) + \partial_\omega
   \varphi^\pm_n(\theta,\omega)\right| \ , \\ \gamma_n^\psi& := & \sup_{\theta
   \in {\cal I}_n+\omega} \left|\partial_\theta
   \psi^\pm_n(\theta,\omega)+\partial_\omega \psi^\pm_n(\theta,\omega)\right|
   \ .
\end{eqnarray}

\begin{lem}
  \label{lem:domgraphs}
Suppose (\ref{eq:Cinvariance})--(\ref{eq:crossing})
  hold and $\omega \in {\cal F}_n$. Then
\begin{equation} \label{eq:domphi}
  \gamma_n^\varphi \ \leq \ S\cdot
 \sum_{k=1}^\infty k \alpha_-^{k}  
\end{equation}
and 
\begin{equation} \label{eq:dompsi}
  \gamma_n^\psi \ \leq \ S \cdot
   \sum_{k=1}^\infty (k+1)\alpha_+^{-k}  \ .
\end{equation}
\end{lem}
\proof 
For any $k,{\cal L} \in \N$ $(\theta,x) \in \ntorus$ there holds
\begin{eqnarray} \nonumber 
  \partial_\omega f^{k+1}_{\theta-({\cal L}+1)\omega}(x) & = & -({\cal L}+1-k)
  \cdot \partial_\theta f_{\theta-({\cal L}+1-k)\omega}(f^k_{\theta-({\cal
  L}+1)\omega}(x)) \\ \label{eq:domstepI} & + & \partial_x f_{\theta-({\cal
  L}+1-k)\omega}(f^k_{\theta-({\cal L}+1)\omega}(x)) \cdot \partial_\omega
  f^k_{\theta-({\cal L}+1)\omega}(x) \ .
\end{eqnarray}
As in the preceding proof, let $(\theta_0,x_0) = (\theta-M_n\omega,c^\pm)$ and
${\cal L} = M_n-1$. Then (\ref{eq:domstepI}) simplifies to
\begin{equation} \label{eq:domstepII}
 \partial_\omega f^{k+1}_{\theta_0}(x_0) \ = \ -({\cal L}+1-k) \cdot
 \partial_\theta f_{\theta_k}(x_k) + \partial_x f_{\theta_k}(x_k) \cdot
 \partial_\omega f^k_{\theta_0}(x_0) \ ,
\end{equation}
and inductive application gives
\begin{equation} \label{eq:domchainI}
  \partial_\omega f^{{\cal L}+1}_{\theta_0}(x_0) \  = \ -\partial_\theta
  f_{\theta_{\cal L}}(x_{\cal L}) - \sum^{{\cal L}-1}_{k=0}({\cal L}+1-k) \cdot \partial_x
  f^{{\cal L}-k}_{\theta_{k+1}}(x_{k+1}) \cdot \partial_\theta
  f_{\theta_{k}}(x_{k}) \ .
\end{equation}
Combining this with (\ref{eq:dthchainI}) and using (\ref{eq:dxestimate})
yields
\begin{eqnarray}\label{eq:dom-calcI}
\lefteqn{\left|\partial_\theta \varphi^\pm_n(\theta,\omega) + \partial_\omega
   \varphi^\pm_n(\theta,\omega)\right| \ = \ \left|\partial_\theta f^{{\cal
   L}+1}_{\theta_0}(x_0) + \partial_\omega f^{{\cal
   L}+1}_{\theta_0}(x_0)\right| } \\\label{eq:dom-calcII} & = & \left| \sum_{k=0}^{{\cal L}-1}
   ({\cal L}-k) \cdot \partial_x f^{{\cal L}-k}_{\theta_{k+1}}(x_{k+1}) \cdot
   \partial_\theta f_{\theta_k}(x_k) \right| \\ \label{eq:dom-calcIII}& \leq & \sum_{k=0}^{{\cal
   L}-1} ({\cal L}-k) \cdot \alpha_-^{{\cal L}-k} \cdot S \
   \leq \ S\cdot \sum_{k=1}^\infty k \alpha_-^{k} \ .
\end{eqnarray}
This proves (\ref{eq:domphi}). 

\ \\
Now let $(\theta_0,x_0) = (\theta+M_n\omega,e^\pm)$ and ${\cal R} =
M_n$. Similar to (\ref{eq:domchainI}) there holds
\begin{equation} \label{eq:dom-calc-backwards}
  \partial_\omega f_{\theta_0}^{-\cal R}(x_0) \ = \ \sum_{k=0}^{{\cal R}-1}
  ({\cal R}-k) \cdot \partial_x f^{-{\cal R}+k+1}_{\theta_{-k-1}}(x_{-k-1})
  \cdot \partial_\theta f^{-1}_{\theta_{-k}}(x_{-k}) \ .
\end{equation}
Using (\ref{eq:dthdxproductbound}) as in the proof of Lemma~\ref{lem:dthgraphs} we obtain
\begin{eqnarray*} 
  \lefteqn{\left|\partial_\omega \psi_n^\pm(\theta,\omega)\right| \ = \
   |\partial_\omega f_{\theta_0}^{-\cal R}(x_0)| } \\ & \leq &  \sum_{k=0}^{{\cal
   R}-1} ({\cal R}-k) \cdot \alpha_+^{-({\cal R}-k)} \cdot 
   S \ \leq \ S \cdot \sum_{k=1}^\infty
   k\alpha_+^{-k} \ .
\end{eqnarray*}
Combined with (\ref{eq:dthpsi}), this yields (\ref{eq:dompsi}). 

\qed

\begin{figure}[h!]  
\noindent 
\hspace{1eM}
\begin{minipage}[t]{\linewidth}  
  \epsfig{file=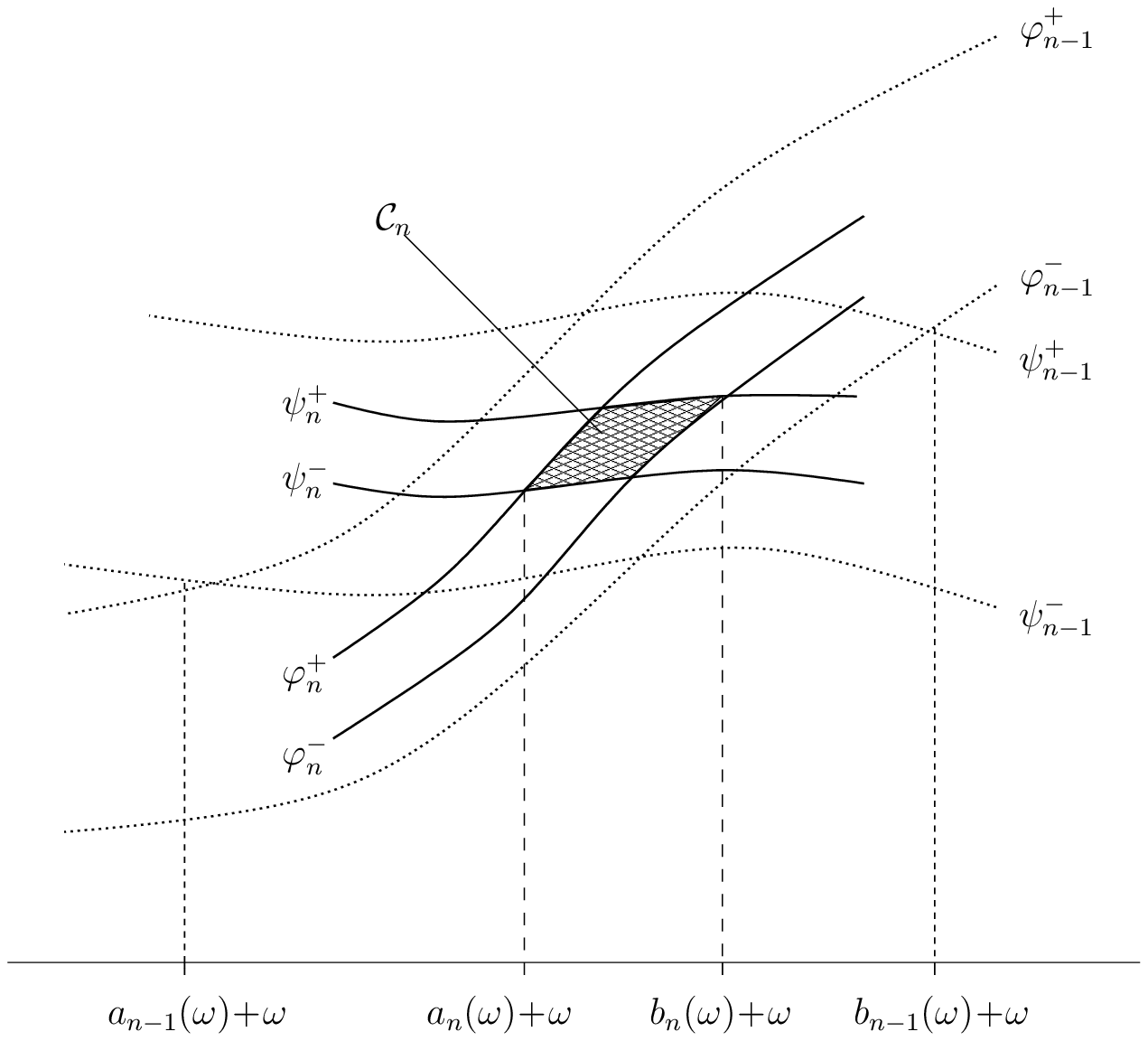, clip=, width=0.85\linewidth}
         \caption{\small The intersection of $f^{M_n}({\cal A}_n^\iota)$ and
         $f^{-M_n}({\cal B}_n^\iota)$.}
         \label{f.intersection}
\end{minipage}
\end{figure}

\begin{lem} \label{lem:Insize}
  Suppose that (\ref{eq:Cinvariance}) holds and $\omega \in {\cal
  F}_n$. Further assume that \ref{eq:In} and \ref{eq:phipsin} hold and
  $l^\varphi_n>u^\psi_n$. Then $({\cal I})_{n+1}$ and $(\Phi/\Psi)_{n+1}$ hold
  and for all $\iota = 1 \ld {\cal N}$. In addition
\begin{equation} \label{eq:In+1length}
  \frac{h^\varphi_n+h^\psi_n}{u^\varphi_n+u^\psi_n} \ \leq \ |I_{n+1}^\iota| \
   \leq \ \frac{H^\varphi_n + H^\psi_n}{l^\varphi_n-u^\psi_n}
\end{equation}
and 
\begin{equation} \label{eq:dwIn+1}
  |\partial_\omega I_{n+1}^\iota| \ \leq \ \frac{\gamma^\varphi_n +
   \gamma_n^\psi}{l^\varphi_n-u^\psi_n} \ .
\end{equation}
\end{lem}
\proof  As $f^{M_n}({\cal A}_n^\iota) \ssq
f^{M_{n-1}}({\cal A}_{n-1}^\iota)$ and $f^{-M_n}({\cal B}_n^\iota) \ssq
f^{-M_{n-1}}({\cal B}_{n-1}^\iota)$ (see Cor.\ \ref{cor:keyestimate}),
\ref{eq:phipsin} implies
\begin{eqnarray*}
J^\varphi_n(a_n(\omega)+\omega) \cap J^\psi_n(a_n(\omega)+\omega)
 & = &
\emptyset \ , \\  J^\varphi_n(b_n(\omega)+\omega) \cap J^\psi_n(b_n(\omega)+\omega)
& = & \emptyset \ .
\end{eqnarray*}
As $|\partial_\theta \varphi^\pm_n - \psi^\pm_n| \geq l^\varphi_n-u^\psi_n >
0$ by assumption, this ensures that the intersection has the geometry depicted
in Figure~\ref{f.intersection}. Hence it is obvious that $I^\iota_n$ contains exactly one
connected component $I^\iota_{n+1}$ of ${\cal I}_{n+1}$, which is not reduced
to a single point. Since ${\cal I}_{n+1}$ is open by definition, this implies
the first two statements of $({\cal I})_{n+1}$. In addition
$I^\iota_{n+1}(\omega) = (a_{n+1}(\omega),b_{n+1}(\omega))$ is characterised
by the equations
\begin{eqnarray*}
     \varphi^+_n(a_{n+1}(\omega)+\omega,\omega) & = &
   \psi_n^-(a_{n+1}(\omega)+\omega,\omega)\ , \\
   \varphi^-_n(b_{n+1}(\omega)+\omega,\omega) & = &
   \psi_n^+(b_{n+1}(\omega)+\omega,\omega) \ ,
\end{eqnarray*}
which yields $(\Phi/\Psi)_{n+1}$.  Further, the estimates
(\ref{eq:phidistance}) and (\ref{eq:psidistance}) in
Lemma~\ref{lem:graphdistance} imply
\[
h^\varphi_n+h^\psi_n \ \leq \
\psi^+_n(a_{n+1}(\omega)+\omega,\omega)-\varphi^-_n(a_{n+1}(\omega)+\omega,\omega)
\ \leq \ H^\varphi_n+H^\psi_n \ ,
\]
and from Lemma~\ref{lem:dthgraphs} we obtain
\[
l^\varphi_n-u^\psi_n \ \leq \ \partial_\theta(\varphi^-_n-\psi^+_n) \ \leq \
u^\varphi_n+u^\psi_n \ . 
\] 
(Note that the bounds in these two lemmas do not depend on $\omega \in {\cal
F}_n$.)  Together, this yields (\ref{eq:In+1length}).

In order to prove (\ref{eq:dwIn+1}), we apply the implicit function theorem to the
identity
\[
 \varphi^+_n(a_{n+1}(\omega)+\omega,\omega) -
 \psi^-_n(a_{n+1}(\omega)+\omega,\omega) \ = \ 0 \ ,
\]
and obtain
\[
\partial_\omega a_{n+1}(\omega) \ = \ \frac{ (\partial_\theta +
   \partial_\omega)\varphi^+_n(a_{n+1}(\omega)+\omega,\omega) - (\partial_\theta +
   \partial_\omega)\psi^-_n(a_{n+1}(\omega)+\omega,\omega)}{\partial_\theta
   \varphi_n^+(a_{n+1}(\omega)+\omega,\omega) - \partial_\theta
   \psi_n^-(a_{n+1}(\omega)+\omega,\omega)}
\]
Therefore (\ref{eq:dwIn+1}) follows from the definitions of
$\gamma^\varphi_n,\gamma^\psi_n,l^\varphi_n$ and $u^\psi_n$, with the same
argument applied to $b_{n+1}$. Consequently ${\cal I}_{n+1}$ depends
differentially on $\omega \in {\cal F}_n$, and the fact that the set ${\cal
F}_{n+1}$ is open follows quite easily from its definition. Thus $({\cal
I})_{n+1}$(iii) holds as well, and this completes the proof.

\qed
\bigskip

We summarise the results of this section in the following proposition, which
is already adapted for its use in the later sections. This is also the reason
why we make the dependence of ${\cal F}_n$ on $M_0 \ld M_n$ explicit in the
statement.

\begin{prop} \label{prop:geometryconclusion}
Suppose (\ref{eq:Cinvariance})--(\ref{eq:crossing}) hold and let $\omega \in {\cal
F}_n(M_0\ld M_n)$. Further, assume that
\begin{equation} \label{eq:calS-cond}
  {\cal S} \ := \ s-S\cdot\left(\frac{1}{\alpha_-^{-1}-1} +
  \frac{1}{\alpha_+-1}\right) \ \geq \ \frac{s}{2}
\end{equation}
and 
\begin{equation}
  \label{eq:gamma-cond}
\gamma \ := \ S \cdot \sum_{k=1}^\infty \left(k\alpha_-^{k} +
(k+1)\alpha_+^{-k}\right) \ \leq \ \frac{\cal S}{4} \ .
\end{equation}
Then $({\cal I})_{n+1}$ and $(\Phi/\Psi)_{n+1}$ hold and for all $j=1 \ld n+1$
and $\iota=1\ld{\cal N}$ we have
\begin{eqnarray}
  |I^\iota_j| & \leq & \frac{2}{s} \cdot
  \max\{\alpha_-,\alpha_+^{-1}\}^{M_{j-1}} \ , \\ |\partial_\omega I^\iota_j|
  & \leq & \viertel \ .
\end{eqnarray}
\end{prop}
\proof Suppose that (\ref{eq:Cinvariance})--(\ref{eq:crossing}) hold. As
already mentioned before, $({\cal I})_0$ and $(\Phi/\Psi)_0$ follow directly
from (\ref{eq:Cinvariance}) and the definition of ${\cal F}_0$. We proceed by
induction.

 Assume that \ref{eq:In} and \ref{eq:phipsin} hold for some $n \geq 0$ and
$\omega \in {\cal F}_n$. Due to Lemma~\ref{lem:dthgraphs} and
(\ref{eq:calS-cond}) we have $l^\varphi_n-u^\psi_n \geq {\cal S} \geq s/2 > 0$
. Therefore we can apply Lemma~\ref{lem:Insize}, which implies that $({\cal
I})_{n+1}$, $(\Phi/\Psi)_{n+1}$ hold. Hence, the required estimates on
$|I^\iota_j|$ and $|\partial_\omega I^\iota_j|$ follow from
Lemma~\ref{lem:Insize}, in combination with Lemma~\ref{lem:graphdistance},
Lemma~\ref{lem:domgraphs} and the estimates provided by (\ref{eq:calS-cond})
and (\ref{eq:gamma-cond}).

\qed


\subsection{Good frequencies} \label{GoodFrequencies}

In order to prove Theorem~\ref{thm:basic}, we will have to show that under the
hypothesis of the theorem there exists a set $\Omega \ssq \kreis$ of positive
measure with the property that for any $\omega \in \Omega$ one can find a
monotonically increasing sequence \nofolge{M_n(\omega)} of positive integers,
such that
\[
\omega \ \in \ \ncap {\cal F}_n(M_0(\omega) \ld
M_n(\omega)) \ .
\]
The problem is that in order to choose the sequences $M_n(\omega)$ inductively
for a sufficiently large set of $\omega$, we will have to make use of the
estimates on the length of the connected components of ${\cal I}_n$ in
Proposition~\ref{prop:geometryconclusion}~. However, these estimates depend in
turn on the choice of the sequence \nofolge{M_n(\omega)}. In order to overcome
this obstacle, we restrict ourselves to choosing the sequences
\nofolge{M_n(\omega)} from the set
\[
{\cal M} \ := \ \left\{\nofolge{M_n} \mid M_n \in [N_n,2N_n) \ \forall n \in
\N_0 \right\} \ ,
\]
where \nofolge{N_n} is a sequence of positive numbers which is fixed {\emph a
  priori} (for simplicity, we do not assume that the $N_n$ are integers). In
  this way we can verify that all required estimates hold, independent of the
  particular choice of \nofolge{M_n(\omega)} in ${\cal M}$.
\medskip

We remark that the results of this section are completely independent of the
preceeding one. In fact, they do not even involve the dynamics of the
system. We only assume that \nofolge{{\cal I}_n} is a family of subsets of
\kreis, such that ${\cal I}_n$ depends on the integers $M_0 \ld M_{n-1}$ and
on $\omega$ (as before, we keep this dependence implicit). While we will make use of
the notation introduced in
Definition~\ref{def:Fn}, we do not use the fact that the sets ${\cal I}_n$ are
defined dynamically as in Definition~\ref{def:In}~. 

As an inductive assumption, we will suppose that for given sequences
\nofolge{K_n} and \nofolge{\eps_n} in Definition~\ref{def:Fn} and a
monotonically increasing sequence sequence \nofolge{N_n} of integers with $N_0
\geq 2$ the following holds:
\begin{equation} \tag{${\cal N}1$} \label{eq:N1}
\begin{array}{l}  \textrm{If } n \in \N_0, \  M_j \in
  [N_j,2N_j) \ \forall j \in [1,n] \textrm{ and } \omega \in {\cal F}_n(M_0
    \ld M_n) \textrm{, then }\\ \\ \begin{array}{cll} (i) & ({\cal I}1)_{n+1}
    \textrm{holds, } \\ \\ (ii) & |I^\iota_j| \leq \eps_j & \forall
    j\in[0,n+1], \iota \in [1,{\cal N}]\\ \\ (iii) & |\partial_\omega
    I^\iota_j| \leq \viertel & \forall j \in [0,n+1], \iota \in [1,{\cal N}]
    \end{array}
    \end{array}
\end{equation}
Finally, we assume that 
\begin{equation} \label{eq:N2}
  N_0  \geq  3 \quad \textrm{ and } \quad  \tag{${\cal N}2$} N_{n+1}  > 
  2K_n N_n \ \forall n \in \N_0 \ .
\end{equation}
\begin{lem} \label{lem:parameterexclI} 
  Suppose (\ref{eq:N1}) and (\ref{eq:N2}) hold and let $M_j\in [N_j,2N_j)$ be
  fixed for $j\in[0,n]$. Further assume that
  \begin{equation} \label{eq:Kjcond1} \tag{${\cal K}$}
    \sum_{j=0}^\infty \frac{1}{K_j} \ < \ \frac{1}{6{\cal N}^2} \ .
  \end{equation}
Then for every $\omega \in {\cal F}_n(M_0 \ld M_n)$ there exists an integer $M
\in [N_{n+1},2N_{n+1})$ such that
\[
     d(({\cal I}_{n+1} - (M-1)\omega) \cup ({\cal I}_{n+1} +
     (M+1)\omega),{\cal Y}_n) \ > \ \eps_n \ .
\]
\end{lem}
\proof
If $j \in [0,n]$ then ${\cal I}_{n+1} \ssq {\cal I}_j$ and $\eps_n \leq
\eps_j$. Therefore
\begin{equation} \label{eq:F2lemma1}
d\left({\cal I}_{n+1} -(p-1)\omega,\bigcup_{k=-M_j+1}^{M_j+1} {\cal I}_j +
k\omega\right) \ \leq \ \eps_n
\end{equation}
implies
\begin{equation} \label{eq:F2lemma2}
d\left({\cal I}_j -(p-1)\omega,\bigcup_{k=-M_j+1}^{M_j+1} {\cal I}_j + k\omega\right) \ \leq
\ \eps_j \ .
\end{equation}
We are going to estimate the number of integers in $(N_{n+1},N_{n+1}+2K_nM_n]$
  $\ssq$ $[N_{n+1},2N_{n+1})$ for which (\ref{eq:F2lemma2}) can happen. Due to
  \ref{eq:F1n} and (\ref{eq:N1})(ii), for any $j\in[0,n],\ \iota,\kappa \in
  [1,{\cal N}]$ and any interval $J \ssq \Z$ of length $|J|\leq 2K_jM_j$,
  there is at most one $p\in J$ such that $ d(
  I^\iota_j-(p-1)\omega,I^\kappa_j) \ \leq \ \eps_j$.  Hence, there are at
  most $2M_j+1$ integers $p$ in $J$ such that
\begin{equation}  \label{eq:F2lemma3}
d\left( I^\iota_j-(p-1)\omega,\bigcup_{k=-M_j+1}^{M_j+1} I^\kappa_j + k\omega \right)
\ \leq \ \eps_j \ ,
\end{equation}
and consequently, due to \ref{eq:In}(i), at most ${\cal N}^2(2M_j+1)$ integers
$p$ in $J$ such that
\begin{equation} \label{eq:F2lemma4}
d\left({\cal I}_j-(p-1)\omega,\bigcup_{k=-M_j+1}^{M_j+1} {\cal I}_j + k\omega \right)
\ \leq \ \eps_j \ .
\end{equation}
Dividing the interval $(N_{n+1},N_{n+1}+2K_nM_n]$ into subintervals of length
  $2K_jM_j$, plus maybe one shorter, we obtain that the number of $p$ in
  $(N_{n+1},N_{n+1}+2K_nM_n]$ for which (\ref{eq:F2lemma2}) holds is bounded by
\[
\left(\frac{K_nM_n}{K_jM_j} + 1\right) {\cal N}^2 (2M_j+1) \ \leq
\  \frac{6K_nM_n{\cal N}^2}{K_j} \ .
\]
Summing up over all $j$, this yields that there are at most 
\[
       2K_nM_n \cdot 3{\cal N}^2 \cdot \sum_{j=0}^n \frac{1}{K_j}
\]
$p$ in $(N_{n+1},N_{n+1}+2K_nM_n]$ with $d({\cal I}_{n+1}-(p-1)\omega,{\cal
Y}_n) \leq \eps_n$. Repeating this argument yields the same bound for the
number of $p$ in $(N_{n+1},N_{n+1}+K_nM_n]$ with $d({\cal
I}_{n+1}+(p+1)\omega,{\cal Y}_n) \leq \eps_n$. Hence, due to
(\ref{eq:Kjcond1}) there must be at least one integer $M \in
(N_{n+1},N_{n+1}+K_nM_n] \ssq (N_{n+1},2N_{n+1}]$ with the required property.

\qed
\bigskip

The following lemma is taken from \cite{bjerkloev:2005}:
\begin{lem} \label{lem:frequenceestimate}
Suppose ${\cal I}={\cal I}(\omega)$ consists of exactly ${\cal N}$ connected
components $I^1 \ld I^{\cal N}$, each of length $|I^\iota| \leq \delta$ and
satisfying $|\partial_\omega I^\iota| \leq \gamma < \halb$. Then for $M \geq
2$ and $\eps > 0$ the set
\[
\left\{ \omega \in \kreis \left| \ d\left({\cal I},\bigcup_{j=1}^M {\cal I}+j\omega
\right) \ < \ \eps \right. \right\}
\]
has measure $\leq 2{\cal N}^2M \frac{\delta + \eps}{1-2\gamma}$ and
consists of at most ${\cal N}^2 M^2-1$ connected components. 
\end{lem}
For any $n\in \N_0$ let
\begin{eqnarray}
  u_{n+1} & := & \ 64 \cdot {\cal N}^2 \cdot K_{n+1} \cdot N_{n+1}^2 \cdot
  \frac{\eps_{n+1}}{\eps_n} \\ v_{n+1} & := & \ \frac{8}{\eps_n} \cdot {\cal
  N}^2 \cdot K_{n+1}^2 \cdot N_{n+1}^3
\end{eqnarray}
Further, let $u_0 := 32 {\cal N}^2K_0N_0\eps_0$ and $v_0 =: 4{\cal N}^2 K_0^2N_0^2$.
\begin{lem} \label{lem:parameterexclusionII} 
 Suppose (\ref{eq:N1}), (\ref{eq:N2}) and (\ref{eq:Kjcond1}) hold and $n \geq
  0$. Let $M_j\in[N_j,2N_j)$ be fixed for $j\in[0,n]$ and assume $\Lambda \ssq
  {\cal F}_n(M_0 \ld M_n)$ is an interval. Then for some $r \leq v_{n+1}$ and
  $\nu = 1 \ld r$ there exist disjoint intervals $\Lambda^\nu \ssq \Lambda$
  and numbers $M^\nu \in [N_{n+1},2N_{n+1})$ such that
\begin{equation}
\Lambda^\nu \ssq {\cal F}_{n+1}(M_0 \ld M_n, M^\nu) \ 
\end{equation}
and
\begin{equation}
   \sum_{\nu=1}^r \Leb(\Lambda^\nu) \ \geq \ \Leb(\Lambda) -   u_{n+1} \ .
\end{equation}
\end{lem}
\proof Obviously $\Lambda$ can be divided into at most
$\frac{2N_{n+1}}{\eps_n}$ intervals $\Gamma^\kappa$ of length $\leq
\frac{2\eps_n}{3N_{n+1}}$. For each $\kappa$, let $\omega^\kappa$ be the
midpoint of $\Gamma^\kappa$. According to Lemma~\ref{lem:parameterexclI},
there exist integers $M^\kappa \in [N_{n+1},2N_{n+1})$, such that
\[
     d(({\cal I}_{n+1} - (M^\kappa-1)\omega^\kappa) \cup ({\cal I}_{n+1} +
     (M^\kappa+1)\omega^\kappa),{\cal Y}_n) \ > \ \eps_n \ .
\]
As $M_j \leq 2N_j < N_{n+1} \ \forall j \in [0,n]$ and $|\partial_\omega I^k_j|
\leq \gamma \leq \viertel \
\forall k,j$ we obtain
\[
     d(({\cal I}_{n+1} - (M^\kappa-1)\omega) \cup ({\cal I}_{n+1} +
     (M^\kappa+1)\omega),{\cal Y}_n) \ > \ 0 \ \ \ \forall \omega \in
     \Gamma^\kappa \ .
\]
Thus $({\cal F}2)_{n+1}$ holds for all $\omega \in \Gamma^\kappa$. 

Let $\tilde{\Gamma}^\kappa$ be the set of those $\omega$'s in $\Gamma^\kappa$
that satisfy $({\cal F}1)_{n+1}$. We have to estimate the size and the number
of connected components of $\tilde{\Gamma}^\kappa$. However, since it follows
from (\ref{eq:N1})(i) and (ii) that ${\cal I}_{n+1}$ consists of ${\cal N}$
connected components of length $\leq \eps_{n+1}$ and $|\partial_\omega
I^\iota_{n+1}| \leq \viertel \ \forall \iota \in [1,{\cal N}]$ by
(\ref{eq:N1})(iii), Lemma \ref{lem:frequenceestimate} with
$\delta=\eps_{n+1}$, $\eps=3\eps_{n+1}$, $\gamma=\viertel$ (see
(\ref{eq:N1})(iii)) and $M=2K_{n+1}N_{n+1}$ yields
\[
     \Leb(\tilde{\Gamma}^\kappa \smin \Gamma^\kappa) \ \leq \ 32{\cal N}^2
     K_{n+1}N_{n+1} \eps_{n+1} \ ,
\]
and the number of connected components of $\Gamma^\kappa$ is at most $4{\cal
  N}^2K_{n+1}^2N_{n+1}^2$. Summing up over $\kappa$ yields the statement. 

\qed
\bigskip

Let 
\[
V_{-1} := 1 \quad \textrm{ and } \quad V_n := \prod_{i=0}^n v_i \quad \forall n
\geq 0 \ .
\] 
\begin{prop}
  Suppose (\ref{eq:N1}), (\ref{eq:N2}) and (\ref{eq:Kjcond1}) hold and
  \begin{equation} 
       \sigma \ := \ 1 - \sum_{n=0}^\infty V_{n-1}u_n  \  .
  \end{equation}
Then there exists a set $\Omega \ssq \kreis$ of measure $\Leb(\Omega) \geq
\sigma$, such that for each $\omega \in \Omega$ there exists a sequence
\nofolge{M_n(\omega)} with the property that
\begin{equation} \label{eq:omegainF}
\omega \ \in \ \bigcap_{n\in\N_0} {\cal F}_n(M_0(\omega) \ld
M_n(\omega))  \ .
\end{equation}
\end{prop}
\proof We are going to construct a nested sequence of sets $\kreis \supseteq \Omega_0 \supseteq
\Omega_1 \supseteq \ldots$ with the following properties:
\romanlist
\item $\Omega_n$ consists of $ \rho_n \leq V_n$ disjoint open intervals
  $\Omega_n^1 \ld \Omega_n^{\rho_n}$.
\item $\Leb(\Omega_n) \ \geq \ 1 - \sum_{i=0}^n V_{i-1} u_i$  
\item For each $i = 1 \ld \rho_n$ there exist numbers $M_0^{n,i} \ld
  M_n^{n,i}$ such that 
  \[
    \Omega_n^i \ssq {\cal F}_n(M_0^{n,i} \ld M_n^{n,i}) \ .
  \]
\item For each $k \leq n$ and each $i \in [1,\rho_n]$ there exists a
  unique $\kappa \in [1,\rho_k]$ such that $\Omega_n^i \ssq
  \Omega_k^\kappa$ and $M_j^{n,i} = M_j^{k,\kappa} \ \forall j=0 \ld k$. 
\listend
For $n=0$ we choose $\Omega_0 = {\cal F}_0$. Recall that this is the set of
all $\omega$ which satisfy condition $({\cal F}1)_0$, and the fact that this
set has all  required properties can be deduced from Lemma~\ref{lem:frequenceestimate}~. 

Now suppose $\Omega_0 \ld \Omega_n$ with the above properties exists. Then for
each $i \in [1,\rho_n]$ we can apply Lemma~\ref{lem:parameterexclusionII} to
the component $\Omega_n^i$ and obtain a union of at most $v_{n+1}$ intervals
with overall measure $\geq m(\Omega_n^i) - u_{n+1}$. Doing this for all the at
most $V_n$ components of $\Omega_n$ yields the required set $\Omega_{n+1}$,
with at most $V_{n+1}= v_{n+1} \cdot V_n$ connected components and measure
$\geq 1-\sum_{i=0}^{n+1} V_{i-1}u_i$.

As the sets $\Omega_n$ form a nested sequence, their intersection $\Omega$ has
measure $\geq \sigma$. Further, for any $\omega \in \Omega$ and $n \in \N$
there exists a unique $i_n \in [1,\rho_n]$ with $\omega \in
\Omega_n^{i_n}$. If we let $M_n(\omega) = M_n^{n,i_n}$, then due to property
(iv) we obtain (\ref{eq:omegainF}).

\qed


\subsection{Proof of Theorem~\ref{thm:basic}, Part A: Existence of SNA} \label{BasicProof}

Suppose that the assumptions of Theorem~\ref{thm:basic} hold. First of all, we
choose the sequence $K_n$ in a way that allows to obtain an lower bound on the
asymptotic expansion and contraction rate, namely
\begin{equation} \label{eq:alphapm-est}
  \min\{\alpha_-^{-1},\alpha_+\} \ \geq \ \alpha^{\ptel} \ .
\end{equation}
In order to do so, we fix $t \in\N$ sufficiently large, such that $t \geq 4$
and
\[
    \frac{2^{-t+2}}{{\cal N}^2} \ \leq \ \log\left(\frac{p^2+2}{p^2+1}\right)
    \ .
\]
Then we let $K_n := 2^{n+t}{\cal N}^2$. Note that this choice satisfies
(\ref{eq:Kjcond1}). We obtain
\[
\beta \ = \ \prod_{n=0}^\infty \left(1-\frac{1}{K_n}\right) \ \geq \
\exp\left(-2\sum_{n=0}^\infty \frac{1}{K_n}\right) \ \geq \
\frac{p^2+1}{p^2+2} \ ,
\]
and this implies 
\[
 \alpha_-^{-1} \ = \ \alpha^{\frac{2}{p}\beta - p(1-\beta)} \ \geq \ \alpha^{\ptel}
 \ .
\]
Similarly we obtain $\alpha_+ \geq \alpha^{\ptel}$, such that
(\ref{eq:alphapm-est}) holds. \bigskip

Now let $N_0 := 3$ and $N_{n+1} := \alpha^{N_n/16p}$. As the sequence $N_n$
grows super-exponentially, (\ref{eq:N2}) holds whenever $\alpha$ is
sufficiently large. Further, let
\[
\eps_0 \ := \ \min_{\iota=1}^{\cal N} |I^\iota_0| \quad \textrm{ and } \quad
\eps_n \ := \ \frac{2}{s}\cdot\alpha^{-N_{n-1}/p}\ .
\]
Again, if $\alpha$ is sufficiently large, then on the one hand $\eps_n \geq
3\eps_{n+1} \ \forall n \in \N_0$ (which is the only requirement on the
sequence \nfolge{\eps_n} in Definition \ref{def:Fn}), and on the other hand
(\ref{eq:calS-cond}) and (\ref{eq:gamma-cond}) hold. Therefore we can apply
Proposition~\ref{prop:geometryconclusion} to see that (\ref{eq:N1}) holds for
the sets ${\cal I}_n$ given by Definition~\ref{def:In}~. This means that all
assumptions of Proposition~\ref{lem:parameterexclusionII} are met, and we
obtain a set $\Omega \ssq \kreis$ of measure
\begin{equation} \label{e.leb-omega}
  \Leb(\Omega) \ \geq \ 1-\sum_{n=0}^\infty V_{n-1}u_n \ ,
\end{equation}
with the property that for all $\omega \in \Omega$ there exists a sequence
$\nfolge{M_n(\omega)}$, such that $\omega \in \bigcap_{n\in\N_0} {\cal
    F}_n(M_0(\omega) \ldots M_n(\omega))$. Proposition~\ref{prop:firststep}
then implies that for all $\omega \in \Omega$ the system 
\[
f\thx \ = \ (\thom,f_\theta(x))
\]
has a sink-source-orbit, and consequently a SNA and SNR by
Proposition~\ref{prop:sinksourcesna}~. It remains to estimate the size of
$\Omega$, i.e.\ to obtain a lower bound on the right side of
(\ref{e.leb-omega}).
 \bigskip

In all of the following estimates we assume that $\alpha$ is chosen
sufficiently large, such that in particular the sequence $N_n$ grows
sufficiently fast, and indicate the steps in which this fact is used by
placing $(\alpha)$ over the respective inequality signs. For any $n\in\N_0$ we
have
\begin{eqnarray*}
u_{n+1} & = & 64{\cal N}^2\cdot K_{n+1}\cdot N_{n+1}^2 \cdot
\frac{\eps_{n+1}}{\eps_n} \\ & = & \ 64{\cal N}^2\cdot K_{n+1} \cdot \alpha^{N_n/8p
- N_n/p + N_{n-1}/p} \ \stackrel{(\alpha)}{\leq} \ \alpha^{-3N_n/4p}
\end{eqnarray*}
and
\begin{eqnarray*}
  v_{n+1} & = & \frac{8}{\eps_n} \cdot {\cal N}^2 \cdot K_{n+1}^2 \cdot N_{n+1}^3 \\
  & \leq & 4s \cdot  {\cal N}^2 \cdot K_{n+1}^2 \cdot \alpha^{N_{n-1}/p+3N_n/16p} \
  \stackrel{(\alpha)}{\leq} \ \alpha^{N_n/4p} \ .
\end{eqnarray*}
Now note that 
\[
V_0 \ =  \ v_0 \ = \ {\cal N}^2\cdot K_0^2\cdot N_0^2 \ \stackrel{(\alpha)}{\leq}  \
\alpha^{N_0/4p} \ .
\]
Further, if we suppose that 
\begin{equation} \label{eq:Vn-est}
  V_{n} \ \leq \ \alpha^{N_n/4p}
\end{equation}
then 
\[
V_{n+1} \ = \ V_n \cdot v_{n+1} \ \leq \ \alpha^{N_n/4p+N_n/4p} \
\stackrel{(\alpha)}{\leq}  \ \alpha^{N_{n+1}/4p} \ .
\]
Consequently, by induction, (\ref{eq:Vn-est}) holds for all $n \geq 1$. We
conclude
\[
V_n u_{n+1} \ \leq \ \alpha^{-N_n/4p}  \ 
\]
and 
\begin{equation}
 \label{e.omega-size}
1  - \sum_{n=0}^\infty V_{n-1}u_n \ \geq \ 1 - u_0 - \sum_{n=0}^\infty
\alpha^{-N_n/4p}  \ .
\end{equation}
As $u_0 = 32{\cal N}^2K_0N_0\eps_0 \ra 0$ if $\eps_0 \ra 0$, the right side is
arbitrarily close to 1 if $\alpha$ is large and $\eps_0$ is small.

To summarise, this means that we can choose constants $\tilde\alpha_0$ and
$\tilde c_0$ in such a way that all the assumptions on $\alpha$ used above
hold and (\ref{e.omega-size}) is larger than $1-\delta$ whenever $\alpha \geq
\tilde\alpha_0$ and $\eps_0 \leq \tilde{c}_0$. Then $\Leb(\Omega) \geq
1-\delta$, as required. This proves Theorem~\ref{thm:basic}, except for the
minimality.


\subsection{Proof of Theorem~\ref{thm:basic}, Part B: Minimality} \label{Minimality}

We choose $\tilde\alpha_0$ and $\tilde c_0$ as at the end of the preceding section and
suppose $\alpha \geq \tilde\alpha_0$ and $\eps_0 \leq \tilde c_0$. Further, we fix $\omega
\in \Omega$ and the corresponding sequence $\nfolge{M_n} =
\nfolge{M_n(\omega)}$ and let $f\thx = (\thom,f_\theta(x))$ as before. Recall
that $M_n \in [N_n,2N_n)$ and $N_{n+1} = \alpha^{N_n/16p}$.

We start with some preliminary remarks and estimates.  Let $\eps_n$ and $N_n$
be chosen as in the last section. Since $\alpha \geq \tilde\alpha_0$ and $\omega \in
{\cal F}_n \ \forall n \in \N_0$, the assumptions of
Proposition~~\ref{prop:geometryconclusion} are satisfied for all
$n\in\N_0$. Consequently, for all $n\geq 0$ the statements \ref{eq:In} and
\ref{eq:phipsin} hold and
\[
|I^\iota_n| \ \leq \ \eps_n \quad \forall \iota \in [1,{\cal N}] \  .
\]
Let 
\[
 \Theta \ := \ \kreis \smin \bigcup_{n\in\N_0} {\cal Z}_n \ .
\]
Then
\begin{eqnarray*}
\Leb(\Theta) & \geq & 1-\sum_{n=0}^\infty m({\cal Z}_n) \ \geq \ 1 -
\sum_{n=0}^\infty 4N_n \eps_n \\ &&
1-N_0\eps_0-\sum_{n=1}^\infty \frac{8}{s} \cdot \alpha^{N_{n-1}/16p -
  N_{n-1}/p} \ . 
\end{eqnarray*}
We now choose the constants $\alpha_0\geq \tilde\alpha_0$ and $c_0 \leq \tilde
c_0$, such that for
all $\alpha\geq\alpha_0$ and $\eps_0 \leq c_0$ there holds 
\begin{equation}
  \label{e.alpha0c0}
\Leb(\Theta) \ > \ 1-\frac{1}{4(1+p^2)} \ .
\end{equation}

Let 
\[
S^* \ := \ S + \frac{S}{\alpha^{-1}_--1} \ 
\]
and choose a constant $\Lambda > 1$ with the following property:
\begin{quote}
  \em If $\Gamma \ssq \ntorus$ is the graph of a differentiable curve
  $\gamma:I\ra\kreis$, defined on an interval $I\ssq\kreis$, and $\Gamma$ has
  slope at most $S^*$, then $f^n(\Gamma)$ has slope at most $S^* \cdot
  \Lambda^n$. 
\end{quote}
Further, due to the lower bound in (\ref{eq:In+1length}) and the estimates
provided by combining Lemmas~\ref{lem:graphdistance}, \ref{lem:dthgraphs} and
\ref{lem:Insize}, there exist constants $B>0$ and $\lambda>0$, such that for
any $n\geq 0$ and any connected component $I^\iota_n$ of ${\cal I}_n$ there
holds
\[
|I^\iota_n| \ \geq \ B \cdot \lambda^{-N_{n-1}} \ .
\]
Let $\delta_n := B \cdot \lambda^{-N_{n-1}}$. Since
\begin{eqnarray*}
\Leb\left(\bigcup_{k=n+1}^\infty {\cal V}_k \smin {\cal V}_n \right) & \leq &
\sum_{k=n+1}^\infty (M_k+1) \cdot \eps_k \ \leq \ \sum_{k=n+1}^\infty
\frac{4}{s} \cdot \alpha^{N_{k-1}/8p-N_{k-1}/p} \ ,
\end{eqnarray*}
and due to the super-exponential growth of the sequence $N_n$, there exists
$n_0 \geq 0$ such that for any connected component
$I^\iota_n$ of ${\cal I}_n$ there holds
\begin{equation} \label{e.no}
\Leb\left(\bigcup_{k=n+1}^\infty {\cal
V}_k\smin {\cal V}_n \right) \ < \ \delta_n/2\Lambda^{M_{n-1}} \quad \forall n
\geq n_0 .
\end{equation}
By slightly reducing the set $\Theta$ if necessary, it is therefore possible
to find a set $\Theta^* \ssq \Theta$ with the following properties:
\begin{description}
\item[$(\Theta^*1)$]\quad \em $\Leb(\Theta^*) \ > \ 1-\frac{1}{2(1+p^2)}$;
\item[$(\Theta^*2)$] \quad For any $\theta \in \Theta^*$, any $n\geq n_0$ and any
   $\iota\in[1,{\cal N}]$, the forward orbit $\{\theta+n\omega \mid n\geq 0\}$
   is $\delta_n/\Lambda^{M_{n-1}}$-dense in $(I^\iota_n-(M_n-1)\omega) \smin
   \left(\bigcup_{k=n+1}^\infty {\cal V}_k \smin {\cal V}_n\right)$.
\end{description}
Now we come to the key point of the proof. The crucial observation is the fact
that there is a large set of points with dense orbit - minimality will then
follow by rather general arguments. More precisely, we prove the following:
\begin{claim}
  \label{claim.denseorbits}
Suppose $\theta_0 \in \Theta^*\cap(\Theta^*-\omega)$ and $x_0 \in E^c$. Then the
forward orbit of $(\theta_0,x_0)$ is dense in \ntorus.
\end{claim}
\proof\ For any point $(\theta_0,x_0)\in \ntorus$, denote its forward orbit by
  ${\cal O}^+(\theta_0,x_0) := \{(\theta_k,x_k) \mid k\geq 0\}$.  Suppose
  $\theta \in \Theta^*\cap(\Theta^*-\omega)$ and $x\in E^c$. Since $\Theta^*
  \ssq \Theta \ssq {\cal Z}_0^c \ssq {\cal I}_0^c$, we can use
  (\ref{eq:Cinvariance}) to see that $f_\theta(x) \in C$. Therefore, it
  suffices to show that the forward orbit of any point $(\theta_0,x_0)$ with
  $\theta_0 \in \Theta^*$ and $x_0 \in C$ is dense. Fix such $\theta_0$ and
  $x_0$ and any $\iota \in [1,{\cal N}]$. Further, choose $n_0$ as in
  (\ref{e.no}). We proceed in four steps:

\ \\  \underline{\em Step 1:} \quad {\it If $n\geq n_0$, then $\pi_1({\cal
  O}^+(\theta_0,x_0) \cap {\cal A}_n^\iota)$ is
  $\delta_n/\Lambda^{M_{n-1}}$-dense in $I^\iota_n-(M_n-1)\omega$.}
\medskip

Since $\theta_0 \in \Theta^*$, it is not contained in ${\cal Z}_n$ for any
$n\in\N_0$. Hence, it follows from Lemma~\ref{lem:keyestimateI} that $x_m
\notin C$ implies $\theta_m \in {\cal V}_k$ for some $k\in\N_0$. Now
$I_n^\iota-(M_n-1)\omega$ is disjoint from ${\cal V}_n$ by \ref{eq:F1n} and
\ref{eq:F2n}. Therefore $\theta_m \in I_n^\iota-(M_n-1)\omega$ and $x_m\notin
C$ imply $\theta_k \in \bigcup_{k=n+1}^\infty {\cal V}_k\smin {\cal V}_n$. In
other words, $x_m \in C$ whenever $\theta_m \in (I_n^\iota-(M_n-1)\omega)
\smin \left(\bigcup_{k=n+1}^\infty {\cal V}_k\smin {\cal V}_n\right)$. The
statement follows from property $(\Theta^*2)$ of the set $\Theta^*$.

\ \\  \underline{\em Step 2:} \quad {\it There exists an integer $n_1\geq n_0$,
  such that for all $n\geq n_1$ the set $\pi_2({\cal O}^+(\theta_0,x_0)
  \cap (I^\iota_{n+1} + (M_n+1)\omega)\times \kreis)$ is $2^{-n}$-dense in $E$.}
\medskip

Let $n\geq n_0$. With the notation of
  Section~\ref{CriticalSets}, we have
\[
f^{M_{n+1}}({\cal A}^\iota_{n+1}) = \{\thx \mid \theta \in
  I^\iota_{n+1}+\omega,\ x \in
  [\varphi^-_{n+1}(\theta),\varphi^+_{n+1}(\theta)]\} \ .
\]
Due to the estimates (\ref{eq:phidistance}) in Lemma~\ref{lem:graphdistance}
and (\ref{eq:dthphi}) in Lemma~\ref{lem:dthgraphs}, this set is a small strip%
\footnote{By `strip', we just mean a set which is the region between two
  continuous curves, defined on a subinterval of $\kreis$. By the slope of a
  strip we mean the slope (or derivative) of its boundary curves.}
of vertical size at most $\alpha^{-M_{n+1}/p}$ and slope at most $S^*$. As
described in the proof of Lemma~\ref{lem:Insize}, this strip crosses the strip
$f^{-M_n}({\cal B}^\iota_n)$ from below to above (where we assume again that
the crossing is upwards), see Figure~\ref{f.intersection}. This implies that
the strip $A:= f^{M_{n+1}+M_n}({\cal A}^\iota_{n+1})$ crosses the horizontal
strip $B^\iota_n=(I^\iota_n+(M_n+1)\omega)\times E$ in the same way.

From (\ref{eq:bounds1}) and $\alpha_u=\alpha^p$, it follows that $A$ has
  vertical size at most $\alpha^{-M_{n+1}/p+2M_np}$.  Further, it has slope at
  most $S^* \cdot \Lambda^{M_n}$. Since $\pi_1({\cal O}^+(\theta_0,x_0) \cap
  {\cal A}_{n+1}^\iota)$ is $\delta_{n+1}/\Lambda^{M_n}$-dense in
  $I^\iota_{n+1}-(M_{n+1}-1)\omega$ by Step 1, it follows that $\pi_2(A)$ is
  $d_n$-dense in $E$, where
  \[
  d_n = S^*\cdot \delta_{n+1} + \alpha^{-M_{n+1}/p+2M_np} \ .
  \]
  Given the super-exponential growth of the sequence $N_n$ and $M_n$, there exists $n_1
  \geq n_0$, such that $d_n \leq 2^{-n} \ \forall n \geq n_1$. This completes
  Step 2.

 \ \\  \underline{\em
  Step 3:} \quad {\it $\mbox{cl}({\cal O}^+(\theta_0,x_0))$ contains a
  vertical segment $\{\zeta\} \times E$ for some $\zeta \in \Theta-\omega$.}
\medskip

Due to compactness and since the size of the intervals $I_n^\iota$ goes to
zero as $n$ goes to infinity, there exists a strictly increasing sequence
$\ifolge{n_i}$ of integers and a point $\zeta\in \kreis$, such that the
intervals $I_{n_i+1}^\iota + (M_{n_i}+1)\omega$ converge to $\{\zeta\}$ in
Hausdorff distance. It follows from Step~2 that $\{\zeta\} \times E \ssq
\mbox{cl}({\cal O}^+(\theta_0,x_0))$.

\ref{eq:F1n} and \ref{eq:F2n} imply that $I_{n}^\iota + (M_{n}+1)\omega$ is
contained in ${\cal Z}_n^c-\omega$ for all $n\in\N_0$. Since the sets ${\cal
Z}_n^c-\omega$ form a nested sequence of compact sets, it follows that $\zeta$
is contained in $\bigcap_{i=0}^\infty {\cal Z}_{n_i}^c-\omega =
\Theta-\omega$.

\ \\ \underline{\em Step 4:} \quad {\it ${\cal O}^+(\theta_0,x_0)$ is dense in
\ntorus.}
\medskip

Let $x^\pm := f_\zeta(e^\pm)$. Since $\Theta-\omega$ is disjoint from ${\cal
  I}_0 \ssq {\cal Z}_0-\omega$, (\ref{eq:Cinvariance}) implies $x^\pm \in
  C$. Consequently, \ref{eq:B1n} holds for all $(\zeta+\omega,x)$ with $x
  \in [x^+,x^-]$ and all $n\in\N_0$. 

Let ${\cal L}$ be the smallest positive integer such that $\zeta+({\cal L}+1)\omega \in {\cal
  I}_n$. Then we can use \ref{eq:C3n} together with (\ref{eq:bounds1}) and
  (\ref{eq:bounds3}) to conclude that 
\[
  \partial_x f^{\cal L}_{\zeta+\omega}(x) \ \leq \ \alpha_-^{\cal L} \quad
  \forall x \in [x^+,x^-] \ .
\]
It follows that $f^{\cal L}(\{\zeta+\omega\} \times [x^+,x^-])$ is a vertical
segment of size smaller that $\alpha_-^{\cal L}$. Since ${\cal L} \geq M_n-1$
(due to $\zeta+\omega \notin {\cal Z}_n$) and $n$ was arbitrary, this means
that the length of the corresponding iterates of $\{\zeta+\omega\} \times
[x^+,x^-]$ goes to zero as $n$ goes to infinity. Therefore, the orbit of the
segment $\{\zeta+\omega\} \times [x^-,x^+] = f(\{\zeta\}\times E)$ is dense in
\ntorus. Since $\{\zeta\}\times E \ssq \mbox{cl}({\cal O}^+(\theta_0,x_0))$ by
Step 3°, this completes the proof of Step 4 and the claim.

\qed
\bigskip

The preceeding claim implies in particular that $f$ is topologically
transitive. It follows from Proposition~\ref{p.unique-minimal-set} that there
is a unique minimal set $M$. Obviously, $M$ cannot be a continuous invariant curve with
positive Lyapunov exponent, since the complement of such a curve always
contains at least one further minimal set. It follows from
\cite{sturman/stark:2000} that $f$ must support at least one invariant measure
$\mu$ with non-positive vertical Lyapunov exponent, that is 
\begin{equation} \label{e.neg-lyap}
\lambda(\mu) \ := \ \int_{\kreis} \int_{\kreis} \partial_x \log f_\theta(x) \
d\mu_\theta(x) \ d\theta \ \leq \ 0 \ .
\end{equation}
(Here $\mu_\theta$ are the conditional measures with respect to the
$\sigma$-algebra $\pi^{-1}({\cal B}(\kreis))$.) 

We claim that this is only possible if $M$ intersects $(\Theta^* \cap
(\Theta^*-\omega)) \times E^c$. In order to see this, not that due to
$(\Theta^*1)$, the set $(\Theta^* \cap (\Theta^*-\omega))$ has measure $>
1-1/(1+p^2)$. If $\supp(\mu) \ssq M$ and $M$ is disjoint from $(\Theta \cap
(\Theta-\omega)) \times E^c$, it therefore follows from (\ref{eq:bounds1}) and
(\ref{eq:bounds2}) that 
\[
\lambda(\mu) \ > \ \left(1-\frac{1}{1+p^2}\right) \cdot \log(\alpha^{1/p}) -
\frac{1}{1+p^2} \cdot \log(\alpha^p) \ = \ 0 \ ,
\]
contradicting (\ref{e.neg-lyap}). 

It follows that $M$ intersects $(\Theta^* \cap (\Theta^*-\omega)) \times E^c$,
and since all points from the later set have dense orbits by
Claim~\ref{claim.denseorbits} we obtain $M=\ntorus$. This completes the proof of
Theorem~\ref{thm:basic}~.


\subsection{Proof of Corollary~\ref{cor:firstexample}} \label{FirstExample}

Obviously, we just have to check that the assumptions
 (\ref{eq:Cinvariance})--(\ref{eq:crossing}) of Theorem~\ref{thm:basic} with
 $\alpha_l^{-1}=\alpha_u=\alpha^p$ are satisfied for all large $\alpha$. Here
 $p$ is meant to be the same as in (\ref{e.ap-def}).  In all of the following,
 we assume that $\alpha$ is chosen sufficiently large and just indicate by
 $(\alpha)$ whenever this fact is used. \medskip

 Due to (\ref{e.g'}), there exist $\eps > 0$ and $s>0$, such that
\[
|g'(\theta)| \ > \ s \quad \forall \theta \in g^{-1}(B_\eps(1/2)) \ .
\]
We let ${\cal I}_0 := g^{-1}(B_\eps(1/2))$, such that (\ref{eq:s}) holds by
definition. Note that due to (\ref{e.g}), ${\cal I}_0$ is the disjoint union
of a finite number of open intervals.  In addition, by reducing $\eps$ further
if necessary, we can assume that all connected components have length smaller
than $\eps_0$, where $\eps_0 = \eps_0(\delta,p,s,S,{\cal N})$ from
Theorem~\ref{thm:basic} with $S := \max_{\theta\in\kreis}|g'(\theta)|$. Note
that this choice of $S$ automatically implies (\ref{eq:bounddth}).

Further, we define $e^\pm:=\pm\alpha^{-\frac{2p-1}{2p}}$ and
$c^\pm:=\mp\eps/2$, and let $E=[e^-,e^+]$ and $C=[c^-,c^+]$ as before. Then
for large $\alpha$ we have $h_\alpha(\kreis \smin E) \ssq B_{\eps/2}(1/2)$, since
\[
h_\alpha(e^\pm) \ = \ \pm
\pi\left(\frac{a_p(\alpha^{1/2p})}{2a_p(\alpha/2)}\right) \
\stackrel{\alpha \ra \infty}{\longrightarrow} \ \halb \ .
\]
Consequently $f_\theta(\kreis \smin E) \ssq C \ \forall \theta \notin {\cal
  I}_0$, such that (\ref{eq:Cinvariance}) holds. Similarly, the above choices
  imply that (\ref{eq:crossing}) holds  (provided we take $\eps < \halb$). 
\medskip

For any $\thx \in \ntorus$, there holds 
\begin{eqnarray*}
\partial_x f_\theta(x) & = & h_\alpha'(x) \ \geq \ h_\alpha'(1/2) 
\\ &=& \frac{\alpha \cdot a_p'(\alpha/2)}{2a_p(\alpha/2)} \ = \ \frac{\alpha\cdot
(1+(\alpha/2)^p)^{-1} }{2a_p(\alpha/2)} \ \stackrel{(\alpha)}{\geq} \ \alpha^{-p} \ .
\end{eqnarray*}
Similarly, there holds
\[
\partial_x f_\theta(x) \ = \ h_\alpha'(x) \ \leq \ h_\alpha'(0) \ = \
\frac{\alpha}{a_p(\alpha/2)} \ \stackrel{(\alpha)}{\leq} \ \alpha^p \ .
\]
Thus (\ref{eq:bounds1}) holds. \medskip

Finally, we check~(\ref{eq:bounds2}) and (\ref{eq:bounds3}). Suppose $x\in
E$. Then 
\[
\partial_x f_\theta(x) \ \geq \ h_\alpha'(x) \ = \ \frac{\alpha a_p'(\alpha
  e)}{2a_p(\alpha/2)} \ = \ \frac{\alpha \cdot
  (1+\alpha^{1/2})^{-1}}{2a_p(\alpha/2)} \ \stackrel{(\alpha)}{\geq} \ \alpha^{1/p} \ .
\]
Similarly, if $x\in C$ there holds
\[
\partial_x f_\theta(x) \ \leq \ h_\alpha'(\eps) \ = \ \frac{\alpha \cdot
  (1+(\alpha \eps)^{p})^{-1}}{2a_p(\alpha/2)} \ \stackrel{(\alpha)}{\leq} \
  \alpha^{-1/p} \ .
\]

It follows that for sufficiently large $\alpha$ all the assumptions of
Theorem~\ref{thm:basic} are satisfied. This completes the proof of the
corollary.


\section{Proof of the refined statement for the qpf Arnold circle map} \label{Proofs}

\subsection{Proof of Theorem~\ref{thm:refined}} \label{RefinedTheorem}

In this section, we describe how the basic construction has to be modified in
order to prove Theorem~\ref{thm:refined}~. In fact, only minor changes are
needed. The only thing which has to be done is to improve some of the
estimates in Section~\ref{CriticalSets}, taking advantage of the additional
assumption (\ref{eq:refinedbounddth}), and then adapt the proof from
Section~\ref{BasicProof} accordingly. We remark that all results of
Sections~\ref{Sink-Source-Orbits} and \ref{CriticalSets} only depend on the
assumptions (\ref{eq:Cinvariance})--(\ref{eq:crossing}) and not on the fact
that the parameter $\alpha$ is chosen very large. Therefore, they all
apply in the situation of Theorem~\ref{thm:refined}~. Similarly, we can
still use all results of Section~\ref{GoodFrequencies}, since these were
completely independent of the dynamics.
\medskip

First of all, we slightly modify the definition
of the sets ${\cal F}_n$: We replace condition $({\cal F}1)_0$ by 
\begin{equation} \label{eq:F1'0}
  \tag*{$({\cal F}1')_0$} d\left({\cal I}_0',\bigcup_{k=1}^{2K_0M_0} ({\cal
  I}_0'+k\omega)\right) \ > \ 3\eps_0 \ 
\end{equation}
and define ${\cal F}_n'$ as the set of all frequencies $\omega \in \kreis$
which satisfy \ref{eq:F1'0}, $({\cal F}2)_0$ and $({\cal F}1$-$2)_j \ \forall
j=1\ld n$. Since ${\cal I}_0 \ssq {\cal I}_0'$, condition \ref{eq:F1'0} is
stronger that $({\cal F}1)_0$, which means that ${\cal F}_n' \ssq {\cal
F}_n$. Consequently, all the results from
Sections~\ref{Sink-Source-Orbits}--\ref{GoodFrequencies} remain true if ${\cal
F}_n$ is replaced by ${\cal F}_n'$ in the respective statements.
\medskip

Since the expansion and contraction rates in Theorem~\ref{thm:refined} are
fixed, we have to improve the estimates from Section~\ref{CriticalSets},
making use of the strengthened condition~\ref{eq:F1'0} together with the
additional assumption (\ref{eq:refinedbounddth}). As the proofs are just
slight variations of the corresponding ones in Section~\ref{CriticalSets}, we
keep the exposition rather brief and only describe the needed modifications. First of
all, Lemma~\ref{lem:dthgraphs} will be replaced by the following:

\begin{lem} \label{lem:dthgraphs'}
 Suppose  (\ref{eq:Cinvariance})--(\ref{eq:refinedbounddth})
  hold and $\omega \in {\cal F}_n'$. Then 
\begin{equation}\label{eq:dthphi'}
 s - \frac{s' + \alpha_-^{M_0}S}{\alpha_-^{-1}-1} \ \leq \ l^\varphi_n \ \leq
\ u^\varphi_n \ \leq \ S + \frac{s' + \alpha_-^{M_0}S}{\alpha_-^{-1}-1}
\end{equation}
and 
\begin{equation}
 \label{eq:dthpsi'} 
u^\psi_n \ \leq \ \frac{s' + \alpha_+^{M_0}S}{\alpha_+-1} \ .
\end{equation}
\end{lem}
\proof As in the proof of Lemma~\ref{lem:dthgraphs}, we fix $\theta \in
I^\iota_n+\omega$ and first let $(\theta_0,x_0) = (\theta-M_n\omega,c^\pm)$
and ${\cal L} = M_n-1$, such that $f_{\theta_0}^{{\cal L}+1}(x_0) =
\varphi_n^\pm(\theta)$. We obtain
\begin{eqnarray*}
  \partial_\theta \varphi^\pm_n(\theta) & = & \partial_\theta f^{{\cal
      L}+1}_{\theta_0}(x_0) \\ & \stackrel{(\ref{eq:dthchainI})}{=} &
      \partial_\theta f_{\theta_{\cal L}}(x_{\cal L}) \ + \ \sum_{k=0}^{{\cal
      L}-1} \partial_x f^{{\cal L}-k}_{\theta_{k+1}}(x_{k+1}) \cdot
      \partial_\theta f_{\theta_{k}}(x_{k}) \\ & \stackrel{({\cal
      F}1')_0}{\geq} & s - \sum_{k={\cal L}-M_0}^{{\cal L}-1} \alpha_-^{{\cal
      L}-k} s' - \sum_{k=0}^{{\cal L}-M_0-1} \alpha_-^{{\cal L}-k} S \\ & = &
      s - \frac{s' + \alpha_-^{M_0}S}{\alpha_-^{-1}-1} \ .
\end{eqnarray*}
The second estimate in (\ref{eq:dthphi'}) follows in the same way. \medskip

In order to prove (\ref{eq:dthpsi'}), we can proceed similarly: We let
 $(\theta_0,x_0) = (\theta+M_n\omega,e^\pm)$ and ${\cal R} = M_n$, such that
 $f_{\theta_0}^{-{\cal R}}(x_0) = \psi_n^\pm(\theta)$, and obtain the required
 estimate from (\ref{eq:dthchainII}) and (\ref{eq:dthdxproductbound}) by using
 \ref{eq:F1'0} once more.

\qed
\medskip

Next, we derive an improved version of Lemma~\ref{lem:domgraphs}:
\begin{lem}
  \label{lem:domgraphs'}
Suppose (\ref{eq:Cinvariance})--(\ref{eq:refinedbounddth})
  hold and $\omega \in {\cal F}_n'$. Then
\begin{equation} \label{eq:domphi'}
  \gamma_n^\varphi \ \leq \ s'\cdot
 \sum_{k=1}^\infty k \alpha_-^k +S\cdot \!\!\!\!
 \sum_{k=M_0+1}^\infty k \alpha_-^k 
\end{equation}
and 
\begin{equation} \label{eq:dompsi'}
  \gamma_n^\psi \ \leq \ s' \cdot
   \sum_{k=1}^\infty (k+1)\alpha_+^{-k} + S \cdot \!\!\!\!
   \sum_{k=M_0+1}^\infty (k+1)\alpha_+^{-k}   \ .
\end{equation}
\end{lem}
\proof The proof is almost identical to that of
Lemma~\ref{lem:domgraphs}~. For proving the upper bound on $\gamma^\varphi_n$,
the only difference is that (\ref{eq:refinedbounddth}) is used instead of
(\ref{eq:bounddth}) in order to estimate $|\partial_\theta f_{\theta_k}(x_k)|$
in the last $M_0$ terms of the sum in (\ref{eq:dom-calcII}).

Similarly, the improved bound on $\gamma_n^\psi$ is obtained by using
(\ref{eq:refinedbounddth}) instead of (\ref{eq:bounddth}) when the last $M_0$
terms of the sum on the right side of (\ref{eq:dom-calc-backwards}) are
estimated via (\ref{eq:dthdxproductbound}).

\qed\bigskip

Lemma~\ref{lem:Insize} can be used without any
modifications. Consequently, we arrive at the following conclusion, whose
proof is identical to that of Proposition~\ref{prop:geometryconclusion}~.

\begin{prop} \label{prop:geometryconclusion'}
Suppose (\ref{eq:Cinvariance})--(\ref{eq:refinedbounddth}) hold and let $\omega \in {\cal
F}_n'(M_0\ld M_n)$. Further, assume that
\begin{equation} \label{eq:calS-cond'}
  {\cal S}' \ := \ s-\left(\frac{s'+\alpha_-^{M_0}S}{\alpha_-^{-1}-1} +
  \frac{s'+\alpha_+^{-M_0}S}{\alpha_+-1}\right) \ \geq \ \frac{s}{2}
\end{equation}
and 
\begin{equation}
  \label{eq:gamma-cond'}
\gamma' \ := s'\cdot\sum_{k=1}^\infty \left(k\alpha_-^k +
(k+1)\alpha_+^{-k}\right) + S\cdot \!\!\!\!\sum_{k=M_0+1}^\infty \!\!\!\left(k\alpha_-^k +
(k+1)\alpha_+^{-k}\right) \ \leq \ \frac{{\cal S}'}{4} \ .
\end{equation}
Then $({\cal I})_{n+1}$ and $(\Phi/\Psi)_{n+1}$ hold and for all $j=1 \ld n+1$
and $\iota=1\ld{\cal N}$ we have
\begin{eqnarray}
  |I^\iota_j| & \leq & \frac{2}{s} \cdot \max\{\alpha_-,\alpha_+^{-1}\}^{M_{j-1}} \ , \\
  |\partial_\omega I^\iota_j| & \leq & \viertel \ .
\end{eqnarray}
\end{prop}
\medskip

In order to complete the proof of Theorem~\ref{thm:refined}, we now choose the
sequence $(K_n)_{n\in\N_0}$ as in the proof of Theorem~\ref{thm:basic}, such
that $\alpha_-^{-1},\alpha_+ \geq \alpha^{1/p}$. Further, we let $N_0$ be the
smallest integer larger than $d^{1/4}$. In all of the following, we assume
that $d$ is chosen sufficiently large to ensure all the required
estimates. As before, we define the sequence $(N_n)_{n\in\N}$ recursively by
$N_{n+1}=\alpha^{N_n/16p}$ and let
\[
\eps_0 \ := \ \min_{\iota=1}^{\cal N} |I^\iota_0| \quad \textrm{ and } \quad
\eps_n \ := \ \frac{2}{s}\cdot\alpha^{-N_{n-1}/p} \ .
\] 
If $d_0$ (and consequently $N_0$) is chosen large enough, then
(\ref{eq:N2}) holds and $\eps_n \geq 3\eps_{n+1} \ \forall n\in\N$. Further,
(\ref{eq:calS-cond'}) and (\ref{eq:gamma-cond'}) hold if $d_0$ is large and
$s'/s$ is small (note that the product $\alpha^{-M_0} S \leq \alpha^{-N_0}S$
decays super-exponentially as $d$ is increased). Thus (\ref{eq:N1}) holds by
Proposition~\ref{prop:geometryconclusion'}~. Therefore, we can apply
Proposition~\ref{lem:parameterexclusionII} and obtain
\[
\Leb(\Omega) \ \geq \ 1-\sum_{n=0}^\infty V_{n-1}u_n \ .
\]
From now on the proof is identical to the one of Theorem~\ref{thm:basic}, with
the only difference that the largeness condition on $\alpha$ is replaced by a
largeness condition on $d$ (and thus $N_0$) in all the respective
estimates. In this way, we obtain
\[
\Leb(\Omega) \ \geq \ 1-u_0-\sum_{n=0}^\infty \alpha^{-N_n/4p} \ .
\]
If $d$ goes to infinity, then due to (\ref{e.eps<1/Ad}) and the choice of
$N_0$ the right side tends to 1 (recall that $u_0 = 32{\cal N}^2
K_0N_0\eps_0$). 
\smallskip

The proof of minimality given in Section~\ref{Minimality} literally stays the
same. The only thing which has to be noted is that the estimate in
(\ref{e.alpha0c0}) also holds for fixed $\alpha$, provided $N_0\approx d^{1/4}$
is chosen sufficiently large.
\smallskip

Hence, we can find constants $c_0$ and $d_0$ with the
required property, which completes the proof. 
\medskip

\subsection{Proof of Corollary~\ref{cor:refined}} \label{Proof-Cor-Refined} 

We place ourselves under the hypothesis of the corollary and let 
\[
f_\theta(x) \ := \ h(x) + \beta g_d(\theta) \ ,
\]
where $g_d(\theta) = \cos(2\pi\theta)^d$. Let $C$ and $E$ be chosen as in
(\ref{e.h-cond1}) and (\ref{e.h-cond2}).  First of all, we fix some $\alpha>1$
and choose $p\in \N$ such that $\sup_{x\in C} h'(x) \leq \alpha^{-2/p}$,
$\inf_{x\in E} h'(x) > \alpha^{2/p}$, and in addition $h'(x) \in
(\alpha^{-p},\alpha^p) \ \forall x\in\kreis$. Then $f$ satisfies
(\ref{eq:bounds1})--(\ref{eq:bounds3}).

Let $\eps := \halb d(h(\kreis \smin E),\kreis \smin C)$ and suppose $\beta \in
[1-\eps,1+\eps]$. Define
\[
  {\cal I}_0 \ := \ g_d^{-1}([-1+\eps,-\eps]\cup[\eps,1-\eps]) \ .
\]
Then it is easy to see that $(f_\theta)_{\theta\in\kreis}$ satisfies
(\ref{eq:Cinvariance}) and (\ref{eq:crossing}). Further, since
$$|\partial_\theta f_\theta(x)| \ = \ |\beta g_d'(\theta)| \ = \ | 2\pi\beta d
\cdot \cos(2\pi\theta)^{d-1} \cdot \sin(2\pi\theta)| \ < \ 4\pi d \ ,$$ we can
choose $S$ in (\ref{eq:bounddth}) smaller than $4\pi d$.
\medskip

 We check that $s$ in (\ref{eq:s}) can be chosen in accordance with
(\ref{e.s>d/A}). In order to obtain an estimate $g_d'$ on ${\cal I}_0$, we
check the endpoints of the connected components and the points where
$g_d''(\theta)=0$. Due to the symmetry of $g_d$, we can restrict to the
interval $[0,1/4]$.

First, assume that $g_d(\theta)=\varepsilon$. Then $\cos(2\pi\theta)
=\varepsilon^{1/d}$ and thus $\sin(2\pi\theta)
=\sqrt{1-\varepsilon^{2/d}}$. Hence
$$ g_d'(\theta) \ = \ -2\pi\beta d\cdot
\varepsilon^{(d-1)/d}\sqrt{1-\varepsilon^{2/d}} \ .
$$
Since $a^y=1+\ln(a)y+O(y^2)$ we have
$$ \sqrt{d}\sqrt{1-\varepsilon^{2/d}}\ = \ \sqrt{2\ln\varepsilon+O(1/d)} \ ,
$$
such that for sufficiently large $d$ there holds
$$
|g_d'(\theta)| \ > \ \eps \cdot \sqrt{\ln(\eps)} \cdot \sqrt{d} \ .
$$

Secondly, assume that $g_d(\theta)=1-\varepsilon$. Then
$\cos(2\pi\theta)=(1-\varepsilon)^{1/d}$ and $\sin (2\pi\theta) =
\sqrt{1-(1-\varepsilon)^{2/d}}$. Thus
$$
g_d'(\theta)=-2\pi\beta d\cdot (1-\varepsilon)^{(d-1)/d}\sqrt{1-(1-\varepsilon)^{2/d}}.
$$
Similar as above we conclude that for sufficiently large $d$ there holds
$$
|g_d'(\theta)| \ > \ (1-\eps) \cdot \sqrt{\ln(1-\eps)} \cdot \sqrt{d} \ .
$$

Thirdly, assume that $g_d''(\theta)=0$.  In this case $\sin(2\pi\theta)^2=1/d$
and $\cos(2\pi\theta)^2 = (d-1)/d$. Therefore
$$
g_d'(\theta)=-2\pi \beta d \left(\frac{d-1}{d} \right)^{(d-1)/2}\frac{1}{\sqrt{d}}=
-2\pi\beta\sqrt{d}\left(1-\frac{1}{d} \right)^{(d-1)/2} \ ,
$$
and the last factor is bounded for all $d$.

From the above analysis we conclude that there is a constant $A$, depending
only on $\eps$, such that for all sufficiently large $d$ there holds
$g_d'(\theta) > \sqrt{d}/A$ for all $\theta \in {\cal I}_0$. 
\bigskip

Finally, we let ${\cal I}_0' := B_{\frac{1}{\sqrt[3]{d}}}(0) \cup
B_{\frac{1}{\sqrt[3]{d}}}(\halb)$. Since $\cos(2\pi\theta) \leq 1-|\theta|^2$
in a neighbourhood of 0, we obtain that for any $\theta \in [0,\viertel]\smin
{\cal I}_0'$ there holds
$$
|g_d(\theta)| \ \leq \ \left(1-d^{-2/3}\right)^d \
\stackrel{d\to\infty}{\longrightarrow} \ 0
\ .
$$ By symmetry, the same estimate holds on all of $\kreis \smin {\cal I}_0'$.
Therefore ${\cal I}_0 \ssq {\cal I}_0'$ for large $d$. 

Similarly, we obtain that for any $\theta \in \kreis \smin {\cal I}_0'$ there
holds 
$$|g_d'(\theta)| \ \leq \ 2\pi\beta d \left(1-d^{-2/3}\right)^{d-1} 
\ \stackrel{d\to\infty}{\longrightarrow} \ 0 \ . $$
Consequently, we can choose $s'$ in (\ref{eq:refinedbounddth}) as a fixed constant, independent of $d$,
which implies that $s'/s$ converges to $0$ as $d$ is increased. 
\medskip

This shows that for sufficiently large $d$ all assumptions of
Theorem~~\ref{thm:refined} are satisfied, which completes the proof of the
corollary.


\end{document}
%

%
%

%

\end{document}